\documentclass[3p]{elsarticle}

\usepackage[utf8]{inputenc}
\usepackage[english]{babel}

\usepackage{amsmath}
\usepackage{amsthm}
\usepackage{amssymb}
\usepackage{bm}
\usepackage{textcomp}
\usepackage{accents}
\usepackage{color}
\usepackage{nicematrix}
\usepackage{booktabs}  
\usepackage[pdftex,plainpages=false,colorlinks=true,citecolor=RoyalBlue,linkcolor=Red,urlcolor=black,filecolor=black,bookmarksopen=true]{hyperref}
\usepackage{subcaption}

\newdefinition{rmk}{Remark}
\newproof{pf}{Proof}

\numberwithin{equation}{section}

\graphicspath{{figs/}}

\begin{document}

\begin{frontmatter}
	
\title{Stability analysis of the steady solution for incompressible lid-driven cavity flow\tnoteref{label1}}
\tnotetext[label1]{The work was supported by the Russian Science Foundation (grant No. 23-71-30013).}

\author{D.V.~Lomasov\fnref{lab1}}
\ead{lomasovdv@my.msu.ru}

\author{P.N. Vabishchevich\fnref{lab1,lab2}}
\ead{vab@cs.msu.ru}

\address[lab1]{Lomonosov Moscow State University, 1, building 52, Leninskie Gory, 119991 Moscow, Russia}

\address[lab2]{North-Eastern Federal University, 58, ul. Belinskogo, 677000  Yakutsk,Russia}

\begin{abstract}

Steady solutions at high Reynolds numbers are computed using an iterative method that linearizes the nonlinear convective terms. This linearization reduces the original nonlinear boundary-value problem to a sequence of linear problems, one solved at each iteration, thereby improving the numerical robustness of the algorithm and enabling convergence over a wide range of governing parameters.
The stability of steady incompressible fluid flows is studied within the framework of a linear problem for perturbations. The original steady velocity and pressure fields are subjected to small perturbations, after which the equations are linearized with respect to the perturbations. The resulting linear problem allows one to analyze the evolution of small perturbations and to establish whether the steady solutions under consideration are stable with respect to infinitesimal perturbations.
In this work, the computational implementation is based on spatial discretization by the finite element method on triangular meshes.
The stability of the obtained solutions is analyzed by solving a generalized algebraic eigenvalue problem formed after discretization of the linear problem for perturbations. The presence of eigenvalues with positive real part indicates instability, whereas their absence indicates stability. Calculations were performed on a sequence of refined meshes at various Reynolds numbers for the flow in a square cavity.

\end{abstract}

\begin{keyword}
Navier-Stokes equations \sep incompressible viscous flow \sep lid-driven cavity \sep finite element method, linear stability analysis \sep high Reynolds number

\MSC 76D05 \sep 76M10 \sep 65N25  \sep 76E99 
	
\end{keyword}
	
\end{frontmatter}

\section{Introduction}

Steady flows of an incompressible viscous fluid are of considerable interest both from the viewpoint of the mathematical theory of the Navier--Stokes equations and from the viewpoint of the development and analysis of numerical methods.
The classical division of incompressible viscous fluid flows into laminar, transitional, and turbulent regimes is associated primarily with the role of the Reynolds number $\operatorname{Re}$. This dimensionless parameter characterizes the ratio of inertial and viscous effects.
As the Reynolds number increases, the role of convective transport grows, the nonlinearity of the equations intensifies, regions of large velocity gradients form, and the flow structure becomes more complex. At the same time, the steady solution may lose stability with respect to small unsteady perturbations. Therefore, when studying flows at high $\operatorname{Re}$, it is necessary to distinguish three interrelated but substantially different problems: constructing a steady solution, determining its stability, and studying the unsteady dynamics after loss of stability.

A modern systematic exposition of the theory of hydrodynamic stability is given in the book \cite{SchmidHenningson2001}, which considers problems of spectral stability, transient growth, nonlinear transition mechanisms, and the relationship between linear stability and subsequent turbulent dynamics. In a recent review article \cite{avila2023transition}, transition processes are considered simultaneously from the standpoint of stability of the laminar state, spatiotemporal dynamics, and dynamical systems theory. In numerical modeling of the transition to turbulence, it is necessary to take into account a wide range of spatial and temporal scales through the use of detailed computational meshes \cite{moin1998direct}. 

When considering flows at high $\operatorname{Re}$, it is important to distinguish among stable laminar flow, steady but unstable solutions of the Navier--Stokes equations, transitional periodic or quasi-periodic flow, and fully developed turbulence.
The existence and regularity of solutions of the Navier--Stokes equations in a bounded domain at high Reynolds numbers differ substantially for two-dimensional and three-dimensional problems \cite{korobkov2024steady}. For the two-dimensional steady problem, the theory is significantly more advanced. For example, the existence of a solution has been established for bounded multiply connected two-dimensional domains and axisymmetric three-dimensional domains. However, the corresponding existence theorems assert neither uniqueness nor stability. The set of steady solutions is analyzed in a similar way \cite{foias1977structure}. The main computational aspects of finding multiple solutions are reflected, for example, in the book \cite{gelfgat2018computational}.

In computational fluid dynamics, considerable attention is paid to the problems of solving problems at high Reynolds numbers \cite{john2016finite,bernardi2024mathematics}.
At high $\operatorname{Re}$, numerical difficulties arise not only from the nonlinearity of the equations but also from the transition to a convection-dominated regime. The problem is singularly perturbed, with the small parameter $1/\operatorname{Re}$ multiplying the diffusion operator.
Thin boundary layers and locally large gradients arise, so that standard spatial discretizations can produce nonphysical oscillations of the solution. Moreover, the velocity--pressure formulation leads to a saddle-point problem associated with the incompressibility constraint.
Because of this, various classes of methods are developed and actively used in computational practice: mixed finite elements with inf-sup-stable pairs; SUPG/PSPG and other stabilized FEM; pressure-correction and projection schemes; finite-volume methods; spectral and spectral-element methods. For example, in the review \cite{langtangen2002numerical}, fully implicit, artificial-compressibility, penalty, and pressure/velocity-correction approaches are systematized, and the connection between different formulations is emphasized.
When considering two-dimensional problems, the vorticity--stream function variables are often used. For such problems, a fundamental issue is related to the correct specification of vorticity on a solid boundary, which is considered in detail in the review \cite{napolitano1999review}.
  
For steady problems of incompressible fluid dynamics, another problem is fundamentally important: how to seek the nonlinear solution. Possible approaches to solving such nonlinear problems are based on Picard/Oseen iterations, Newton and inexact Newton methods with continuation in $\operatorname{Re}$ \cite{elman2014finite}.
In our work \cite{lomasov2025numerical}, a simple iterative method was constructed on the basis of linearization of the nonlinear convective term, which converges at high Reynolds numbers and makes it possible to identify several solutions.
Approximate solutions of steady Navier--Stokes boundary-value problems are often obtained by solving the unsteady equations.
In the time-marching approach, an approximate steady-state solution is obtained by integrating the equations over a sufficiently long time.
In pseudo-transient continuation methods, the time term is not directly related to the unsteady Navier--Stokes equation and plays the role of a regularizer of the iterative process \cite{kelley1998convergence}.
For cavity flows, such approaches were implemented in the article \cite{botti2015choice}.

The two-dimensional problem of incompressible fluid flow in a square cavity \cite{shankar2000driven} is the best-known benchmark for Navier--Stokes solvers.
In the key work \cite{ghia1982high}, the authors used the vorticity--stream function formulation and obtained solutions up to $\mathrm{Re}=10 \,000$.
The numerical calculations were refined in the article \cite{botella1998benchmark} with allowance for the asymptotics of singular terms of the solution.
These data are still used for comparing computational algorithms.
Among later studies, we note the work \cite{sahin2003novel-1}, in which a fully implicit velocity-based finite-volume formulation and Newton method are used.
The second part of the work \cite{sahin2003novel-2} is devoted to linear stability analysis at high Reynolds numbers.
The most complete results on finding steady solutions for two-dimensional cavity flows are presented in \cite{Li2025Flow}. The authors provide benchmark-quality steady 2D solutions and study the flow structure up to $\operatorname{Re}=30 \,000$.
Similar data were also obtained in the work \cite{hostos2021solving}.
In such studies, it is important to distinguish numerical experiments that identify stable steady solutions from those that find steady solutions that are already unstable.

Numerical methods for analyzing the stability of steady flows can be conventionally divided into two main classes. The first is based on the direct solution of the unsteady Navier--Stokes equations. An initial perturbation is specified relative to the steady solution, after which its further evolution is analyzed. Decay of the perturbation indicates stability, whereas its growth indicates instability. This approach is a natural way to study nonlinear dynamics and allows one to directly obtain periodic and more complex regimes after loss of stability. Direct numerical simulation of unsteady flows is one of the main tools of modern computational fluid dynamics and is used both for studying the development of instability and for calculating established unsteady regimes \cite{moin1998direct,pope2000turbulent,ferziger2002computational}. General issues of numerical solution of the unsteady Navier--Stokes equations and construction of stable algorithms are considered in detail in \cite{elman2014finite,temam1984navier}.

However, in the immediate vicinity of the critical Reynolds number, the growth or decay rate of perturbations becomes small, so that long-time integration is required. Moreover, direct calculation does not allow one to directly determine the complete spectrum of the linear operator and, as a rule, does not allow one to obtain the unstable steady state. The latter circumstance is of fundamental importance in the study of bifurcations: an unstable steady state may exist as a solution of the nonlinear steady problem but not be observed in direct time integration because of the growth of perturbations. General issues of flow stability, transition to unsteady regimes, and the relationship between linear theory and nonlinear dynamics are systematically presented in the monographs \cite{SchmidHenningson2001,drazin2004hydrodynamic}.

The second class of methods is based on direct spectral analysis of the linearized problem. After spatial discretization, a large generalized eigenvalue problem arises. In practice, it is necessary to compute a small number of eigenvalues that determine the instability of the steady state. Such an approach is a standard tool for studying hydrodynamic stability and allows one to directly study the spectrum of the linearized operator and determine the most dangerous perturbations \cite{SchmidHenningson2001,drazin2004hydrodynamic,schmid2007nonmodal}. When computing the dominant part of the spectrum of large sparse matrices, iterative Krylov subspace methods are used in various variants
\cite{saad2011numerical,watkins2007matrix}.
In particular, the application of Arnoldi-type methods to hydrodynamic stability problems is discussed in detail in \cite{schmid2007nonmodal}. The spectral approach allows one to simultaneously determine the critical Reynolds number, the growth rate of the most dangerous perturbation, and the frequency of the arising oscillations.

Two-dimensional flows in a square cavity are benchmark problems of computational fluid dynamics for analyzing hydrodynamic instability. Numerical studies show that the first two-dimensional instability arises at a Reynolds number of the order of eight thousand and has the character of a Hopf bifurcation.
In the work \cite{bruneau2006cavity}, the authors performed a numerical analysis of the critical region of flows near $\operatorname{Re} = 8\,000$ and confirmed the occurrence of a periodic regime by directly solving the unsteady equations.
A benchmark study of unsteady flow inside a square lid driven cavity was carried out in the work \cite{sengupta2025benchmark}.
Using a spectral method in the article \cite{auteri2002stability}, a critical value was obtained in a narrow interval near $8\,018$. A close value ($8\,032$) was obtained in the work \cite{sahin2003novel-2} using the Arnoldi method. It should be noted that the differences between published values are not of fundamental importance.
The overall picture of loss of stability does not change; only the sensitivity of the results to spatial discretization and to the choice of computational fluid dynamics tools is demonstrated.

For the cavity problem, it is also necessary to distinguish between two-dimensional and three-dimensional stability. The critical Reynolds number of about $8\,000$ refers to the stability of the two-dimensional steady solution with respect to two-dimensional perturbations. If, however, three-dimensional perturbations of the two-dimensional base flow are considered, instability can arise at substantially lower Reynolds numbers. For example, studies \cite{ding1998cavity,albensoeder2001three} showed the presence of three-dimensional unstable modes already at Reynolds numbers less than $1\,000$.

In most previous studies, the main attention was paid to determining the first loss of stability of the steady solution and the corresponding critical Reynolds number associated with the first Hopf bifurcation. Such an analysis does not characterize the subsequent change in the spectrum with a further increase in the Reynolds number.
Despite a significant number of studies of the first loss of stability of the flow in a square cavity, the subsequent evolution of the spectrum of the steady solution with increasing Reynolds number has been studied much less systematically.
In this connection, we first note the work \cite{sengupta2011universal}, in which, the first three critical Reynolds numbers are identified from an analysis of unsteady solutions: $7\,934, 8\,188$, and $9\,934$, respectively.
In the article \cite{osada2011numerical}, on the basis of solving an unsteady problem, the first bifurcation was identified near $\operatorname{Re} = 7\,987$.
A secondary Hopf bifurcation of the developing unsteady regime was also obtained near $9\,575$, and the third one near $10\,975$.

The work \cite{lestandi2018multiple} is devoted to direct numerical simulation with reproduction of the nonlinear spatiotemporal dynamics after the first bifurcation and the arising periodic and quasi-periodic regimes, as well as to separate determination of subsequent critical parameters.
On rather coarse meshes, several transitions with increasing Reynolds number are numerically demonstrated, and not only the first loss of stability.
It is interesting to compare the results of these works on the nonlinear dynamics of perturbations with the corresponding spectral problems of linear analysis of hydrodynamic instability.

In the present work, the spectrum of the linearized operator is directly traced relative to the steady solution branch with a successive increase in the Reynolds number. This allows one to identify various eigenmodes responsible for successive losses of stability and to determine the corresponding critical parameter values in a unified spectral formulation. The calculations are performed for two-dimensional flow in a square cavity up to $\operatorname{Re}=10\,500$, which allows one to study the spectral structure of the steady solution not only in the vicinity of the first bifurcation but also in the subsequent region of its instability.

The article is structured as follows. In the next section, the steady boundary value problem in a bounded domain for the incompressible Navier--Stokes equations is formulated. The study of stability with respect to small perturbations is based on the spectral problem. In Section 3, the benchmark problem of analyzing two-dimensional flows in a square cavity is formulated.
The computational algorithm is based on finite element approximation and solution of a partial eigenvalue problem for the corresponding discrete problem. The results of numerical calculations at various values of the Reynolds number are presented in Section 4.
In the conclusion, a summary and conclusions are given.

\section{Problem formulation}

The problem of analyzing the hydrodynamic instability of steady flows of a viscous incompressible fluid in a bounded domain is considered.
The study is based on solving a spectral problem for small perturbations.

\subsection{Steady problem}

A steady flow of an incompressible fluid is considered in a bounded
domain $\Omega \subset \mathbb{R}^d$, $d = 2$, with a sufficiently smooth
boundary $\partial\Omega = \Gamma$. The velocity field $\bm U(\bm x)$ and
the pressure divided by the density, $P(\bm x)$ satisfy the steady
Navier--Stokes system
\begin{align}
    \operatorname{div}\bigl(\bm U \otimes \bm U\bigr)
    + \operatorname{grad} P
    - \frac{1}{\operatorname{Re}} \operatorname{div}\operatorname{grad}\bm U
    &= \bm 0, \quad \bm x \in \Omega, \label{eq:2-1} \\[2pt]
    \operatorname{div}\bm U &= 0, \quad \bm x \in \Omega .\label{eq:2-2}
\end{align}
The Reynolds number $\operatorname{Re}$ is
\[
    \operatorname{Re} = \frac{\rho \, U_0 \, L}{\mu},
\]
where $\rho$ is density, $\mu$ is dynamic viscosity, and $U_0$ and $L$ are
the characteristic velocity and linear size of the problem. The nonlinearity of system
\eqref{eq:2-1}--\eqref{eq:2-2} is due to the convective term
$\operatorname{div}(\bm U \otimes \bm U)$.

System \eqref{eq:2-1}--\eqref{eq:2-2} is closed by inhomogeneous boundary conditions for velocity, when
\begin{equation}
    \bm U = \bm g(\bm x), \quad \bm x \in \Gamma .
    \label{eq:2-3}
\end{equation}
The given function $\bm g(\bm x)$ satisfies the compatibility condition
\begin{equation}
    \int_\Gamma \bm g \cdot \bm n \, d \bm x = 0 ,
    \label{eq:2-4}
\end{equation}
where $\bm n$ is the unit outward
normal to $\Gamma$. Condition \eqref{eq:2-4} expresses the balance of flow rate
through the boundary of the domain and is necessary for the existence of a solution of problem
\eqref{eq:2-1}--\eqref{eq:2-3}.
The pressure $P$ is determined up to an additive constant; therefore, the normalization condition is imposed
\begin{equation}
    \int_\Omega P \, d \bm x = 0.
    \label{eq:2-5}
\end{equation}
The steady solution of problem \eqref{eq:2-1}--\eqref{eq:2-5} is denoted
below by $(\bm U, P)$ and is assumed to be sufficiently smooth: we assume that
$\bm U \in H^2(\Omega)^d$, $P \in H^1(\Omega) \cap L^0_2(\Omega)$, where
\[
L^0_2(\Omega) = \left\{ q \in L_2(\Omega) : \int_\Omega q\,d \bm x = 0 \right\}.
\]

\subsection{Linearization with respect to the steady solution}

To analyze the stability of the steady solution $(\bm U, P)$ to small
perturbations \cite{SchmidHenningson2001,drazin2004hydrodynamic}, we represent the unsteady flow field as the sum of the base
steady solution and a small perturbation:
\begin{equation}
    \bm u(\bm x, t) = \bm U(\bm x) + \bm u'(\bm x, t), \quad
    p(\bm x, t) = P(\bm x) + p'(\bm x, t),
    \label{eq:2-6}
\end{equation}
where $\bm u'(\bm x, t)$ and $p'(\bm x, t)$ are the velocity and pressure perturbations.
The perturbations satisfy homogeneous boundary conditions
\begin{equation}
    \bm u'(\bm x,t) = \bm 0 
    \quad \bm x \in  \Gamma ,
    \quad t > 0.
    \label{eq:2-7}
\end{equation}

Substitution of representation \eqref{eq:2-6} into the unsteady
Navier--Stokes equations
\begin{align*}
    \frac{\partial \bm u}{\partial t}
    + \operatorname{div}\bigl(\bm u \otimes \bm u\bigr)
    + \operatorname{grad} p
    - \frac{1}{\operatorname{Re}} \operatorname{div}\operatorname{grad}\bm u
    &= \bm 0, \\
    \operatorname{div}\bm u &= 0
\end{align*}
taking into account the steady system \eqref{eq:2-1}--\eqref{eq:2-2}
leads to a system for the perturbations. Using the identity
\[
    \operatorname{div}\bigl(\bm u \otimes \bm u\bigr)
    = \operatorname{div}\bigl(\bm U \otimes \bm U\bigr)
    + \operatorname{div}\bigl(\bm U \otimes \bm u'\bigr)
    + \operatorname{div}\bigl(\bm u' \otimes \bm U\bigr)
    + \operatorname{div}\bigl(\bm u' \otimes \bm u'\bigr)
\]
and discarding second-order terms in the perturbations
$\operatorname{div}(\bm u' \otimes \bm u')$, we obtain the linearized system
\begin{align}
    \frac{\partial \bm u'}{\partial t}
    + \operatorname{div}\bigl(\bm U \otimes \bm u'\bigr)
    + \operatorname{div}\bigl(\bm u' \otimes \bm U\bigr)
    + \operatorname{grad} p'
    - \frac{1}{\operatorname{Re}} \operatorname{div}\operatorname{grad}\bm u'
    &= \bm 0, 
    \quad \bm x \in \Omega ,
    \quad t > 0 , \label{eq:2-8} \\[2pt]
    \operatorname{div}\bm u' &= 0, 
    \quad \bm x \in \Omega ,
    \quad t > 0 , \label{eq:2-9}
\end{align}
which is supplemented by homogeneous boundary conditions \eqref{eq:2-7} and
\begin{equation}
    \int_\Omega p' \, d \bm x = 0.
    \label{eq:2-10}
\end{equation}
Given initial perturbations
\begin{equation}
    \bm u'(\bm x,0) = \bm w(\bm x),
    \quad \bm x \in \Omega ,
    \label{eq:2-11}
\end{equation}
the boundary value problem \eqref{eq:2-7}--\eqref{eq:2-11} describes the evolution of small
perturbations in the linear approximation.

\subsection{Spectral problem}

Linear spectral analysis of hydrodynamic stability \cite{SchmidHenningson2001,drazin2004hydrodynamic} is based on the representation of solutions of the linearized boundary value problem \eqref{eq:2-7}--\eqref{eq:2-11}
in the form of normal modes:
\begin{equation}
    \bm u'(\bm x, t) = \bar{\bm u}(\bm x)\, e^{\lambda t}, \quad
    p'(\bm x, t) = \bar{p}(\bm x)\, e^{\lambda t}.
    \label{eq:2-12}
\end{equation}
Here $\lambda \in \mathbb{C}$ is a complex decrement determining the temporal evolution of the mode, and $\bar{\bm u}(\bm x)$ and $\bar{p}(\bm x)$ are the amplitudes of the velocity and pressure perturbations.
Substitution of
\eqref{eq:2-12} into \eqref{eq:2-7}--\eqref{eq:2-10} gives
the spectral problem: find $\lambda \in \mathbb{C}$, $\bar{\bm u}$, and
$\bar{p}$ from the equations
\begin{align}
    \lambda \bar{\bm u}
    + \operatorname{div}\bigl(\bm U \otimes \bar{\bm u}\bigr)
    + \operatorname{div}\bigl(\bar{\bm u} \otimes \bm U\bigr)
    + \operatorname{grad} \bar{p}
    - \frac{1}{\operatorname{Re}} \operatorname{div}\operatorname{grad}\bar{\bm u}
    &= \bm 0, \quad \bm x \in \Omega, \label{eq:2-13} \\[2pt]
    \operatorname{div}\bar{\bm u} &= 0, \quad \bm x \in \Omega, \label{eq:2-14}
\end{align}
under the constraints
\begin{equation}
    \bar{\bm u} = \bm 0, \quad \bm x \in \Gamma,
    \label{eq:2-15}
\end{equation}
\begin{equation}
    \int_\Omega \bar{p} \, dx = 0.
    \label{eq:2-16}
\end{equation}

The spectral problem \eqref{eq:2-13}--\eqref{eq:2-16} can be written in operator form
\begin{equation}
    \mathcal{L}(\bm U) \, \bar{\bm u} + \mathcal{G} \bar{p}
    = \lambda \, \bar{\bm u}, \quad
    \mathcal{D} \bar{\bm u} = 0,
    \label{eq:2-17}
\end{equation}
with the linear operator $\mathcal{L}(\bm U)$ defined as
\begin{equation}
    \mathcal{L}(\bm U) \, \bar{\bm u}
    = -\operatorname{div}\bigl(\bm U \otimes \bar{\bm u}\bigr)
    - \operatorname{div}\bigl(\bar{\bm u} \otimes \bm U\bigr)
    + \frac{1}{\operatorname{Re}} \operatorname{div}\operatorname{grad}\bar{\bm u},
    \label{eq:2-18}
\end{equation}
where $\mathcal{G} = -\operatorname{grad}$ is the gradient operator with the opposite sign, and
$\mathcal{D} = \operatorname{div}$ is the divergence operator. The operator
$\mathcal{L}(\bm U)$ depends on the base steady flow $\bm U$ and
the Reynolds number $\operatorname{Re}$.

Under inhomogeneous boundary conditions \eqref{eq:2-3}, the base flow $\bm U$ is not trivial, and the operator $\mathcal{L}(\bm U)$ defined according to \eqref{eq:2-18} contains both convective transport by the base flow
$\operatorname{div}(\bm U \otimes \bar{\bm u})$ and deformation of the base
flow $\operatorname{div}(\bar{\bm u} \otimes \bm U)$. Both terms contribute to the non-self-adjointness of the operator $\mathcal{L}(\bm U)$ and,
as a consequence, to the possible existence of complex eigenvalues
$\lambda$.

\subsection{Analysis of hydrodynamic stability}

The spectrum of problem \eqref{eq:2-13}--\eqref{eq:2-16} consists of eigenvalues $\lambda_k = \Re \lambda_k \pm i \Im \lambda_k \in \mathbb{C}$. The real part
$\Re \lambda_k$ determines the growth or decay rate of the corresponding mode:
\begin{itemize}
    \item if $\Re \lambda_k < 0$, then the perturbation decays exponentially in
    time and the mode is stable;
    \item if $\Re \lambda_k > 0$, then the perturbation grows exponentially and
    the mode is unstable;
    \item if $\Re \lambda_k = 0$, then neutral stability takes place.
\end{itemize}
The steady flow $(\bm U, P)$ is linearly stable if all eigenvalues of problem \eqref{eq:2-13}--\eqref{eq:2-16} $\lambda_k(\operatorname{Re})$ lie in the left half-plane:
\begin{equation}
    \max_k \Re \lambda_k (\operatorname{Re}) < 0.
    \label{eq:2-19}
\end{equation}
Loss of stability occurs when at least one eigenvalue crosses the imaginary axis:
\[
    \Re \lambda_k (\operatorname{Re}) = 0.
\]
The corresponding value of the Reynolds number is called critical and is denoted by $\operatorname{Re}_{\mathrm{cr}}$.

The aim of our study is the numerical analysis of changes in the spectral stability of the steady incompressible fluid flow $(\bm U, P)$ (problem \eqref{eq:2-1}--\eqref{eq:2-5}) with a successive increase in the Reynolds number. For a separate value of $\operatorname{Re}$, the steady solution of the Navier--Stokes equations (problem \eqref{eq:2-1}--\eqref{eq:2-5}) is first computed, and the steady branch is continued even after the first loss of stability. Then this steady solution is analyzed for stability with respect to small perturbations. For this, the corresponding generalized eigenvalue problem \eqref{eq:2-17}--\eqref{eq:2-18} is solved. For each $\operatorname{Re}$, the dominant eigenvalues closest to the imaginary axis are determined. On the basis of the trajectories of individual eigenvalues, a picture of the evolution of the spectrum of the steady solution with increasing Reynolds number is constructed.
Numerical experiments were carried out for the two-dimensional problem of fluid flow in a square cavity.
The calculations were performed with an increase in $\operatorname{Re}$ up to $10\,500$.

\section{Computational algorithm}

The computational technology of spectral analysis of hydrodynamic stability was tested on the two-dimensional standard lid-driven cavity problem. We describe the key elements of the computational algorithm used to approximate the steady problem at high Reynolds numbers by the finite element method and to solve the spectral problem for the dominant modes.

\subsection{Two-dimensional problem of flow in a square cavity}

As a test problem, we consider the steady flow of an incompressible
fluid in a cavity with cross-section $\Omega = (0,1) \times (0,1) \ (d =2, \ \bm x = (x_1, x_2))$ with a moving upper
lid \cite{ghia1982high}. Despite the simple geometry and
dependence of the solution on only one parameter---the Reynolds number---this problem
reflects almost all main features of incompressible viscous
fluid flows and is widely used for verification of numerical methods~\cite{botella1998benchmark,erturk2005}.

We solve boundary value problem \eqref{eq:2-1}--\eqref{eq:2-5}, in which
\[
    \bm g(\bm x) = \begin{cases}
    (1, 0), & \bm x \in \Gamma, \quad x_2 = 1, \\
    (0, 0), & \bm x \in \Gamma, \quad x_2 < 1 .
    \end{cases} 
\]
As the Reynolds number increases, a complex vortex
structure forms in the cavity: a central vortex, secondary vortices in the corners, and, at sufficiently large
$\operatorname{Re}$, tertiary vortices~\cite{erturk2005}. It is known that the first
Hopf bifurcation, leading to loss of stability of the steady solution,
occurs at $\operatorname{Re} \approx 8\,000$~\cite{sahin2003novel-2, auteri2002stability}.

\subsection{Numerical solution of the steady problem}

To solve the nonlinear system of steady Navier--Stokes equations, an
iterative method based on linearization of the convective term is used. The performance of this variant of the Picard method is discussed in \cite{lomasov2025numerical} for the approximate solution of the two-dimensional steady problem for incompressible fluid flows in a cavity at high Reynolds numbers.

An essential point is related to the choice of linearizations of the convective term. Instead of the usually applied approach (see, for example, \cite{turek1999efficient}) with linearization at each iteration
$n+1$
\[
    \operatorname{div}\bigl(\bm U^{(n+1)} \otimes \bm U^{(n+1)}\bigr)
    \;\approx\;
    \operatorname{div}\bigl(\bm U^{(n)} \otimes \bm U^{(n+1)}\bigr)
\]
for determining $\bm U^{(n+1)}$, we use
\[
    \operatorname{div}\bigl(\bm U^{(n+1)} \otimes \bm U^{(n+1)}\bigr)
    \;\approx\;
    \operatorname{div}\bigl(\bm U^{(n+1)} \otimes \bm U^{(n)}\bigr) .
\]
In addition, to accelerate convergence and improve the stability of the iterative process,
relaxation is applied.

Applied to the problem under consideration \eqref{eq:2-1}--\eqref{eq:2-5}, we first find $\widetilde{\bm U}^{(n+1)}$ and $\widetilde{P}^{(n+1)}$ from
the equations
\[
\begin{split}
    \operatorname{div}\bigl(\widetilde{\bm U}^{(n+1)} \otimes \bm U^{(n)}\bigr)
    + \operatorname{grad} \widetilde{P}^{(n+1)}
    - \frac{1}{\operatorname{Re}} \operatorname{div}\operatorname{grad}\widetilde{\bm U}^{(n+1)}
    &= \bm 0, \quad \bm x \in \Omega, 
     \\[2pt]
    \operatorname{div}\widetilde{\bm U}^{(n+1)} &= 0, \quad \bm x \in \Omega ,
\end{split}
\]
under the additional conditions
\[
    \widetilde{\bm U}^{(n+1)} = \bm g(\bm x), \quad \bm x \in \Gamma ,
    \quad 
    \int_\Omega \widetilde{P}^{(n+1)} \, d \bm x = 0.
\]
After this, the new approximation is determined:
\[
    \bm U^{(n+1)} = \omega \, \widetilde{\bm U}^{(n+1)}
    + (1 - \omega) \, \bm U^{(n)},
\]
\[
    P^{(n+1)} = \omega \, \widetilde{P}^{(n+1)}
    + (1 - \omega) \, P^{(n)},    
\]
where $\omega > 0$ is the relaxation parameter. The choice of $\omega$ is carried out empirically. This variant of the Picard method with relaxation ensures reliable convergence even for high Reynolds numbers.
Moreover, by choosing the relaxation parameter, it is possible to identify two solutions for flows in a semi-elliptical cavity \cite{lomasov2025numerical}.
In the calculations presented below, we used the value $\omega=0.8$.

\subsection{Variational formulation and finite element method}

We multiply the momentum conservation equation by a vector test function
$\bm v \in V =  H_0^1(\Omega)^2$, and the continuity equation by a scalar test function
$q \in Q = L_2^0(\Omega)$. Integration over the domain $\Omega$ leads us to the variational formulation of problem \eqref{eq:2-1}--\eqref{eq:2-5}: find $(\bm U, P)$ taking into account boundary conditions \eqref{eq:2-3}, such that for all $(\bm v, q) \in V \times Q$
\[
\begin{split}
    \frac{1}{\operatorname{Re}} \int_\Omega \operatorname{grad} \bm U : \operatorname{grad} \bm v \, d \bm x
    + \int_\Omega \operatorname{div}\bigl(\bm U \otimes \bm U\bigr) \cdot \bm v \, d \bm x
    - \int_\Omega P \, \operatorname{div}\bm v \, d\bm x &= 0, \\[2pt]
    \int_\Omega q \, \operatorname{div}\bm U \, d \bm x &= 0 .
\end{split}
\]

For spatial discretization, uniform triangular meshes with $N$ subdivisions in each coordinate direction are used. The main computational data were obtained with $N = 1\,024$.
As the finite element pair, Taylor--Hood elements are chosen \cite{taylor1973}, where velocity is approximated by continuous piecewise quadratic polynomials, and pressure by continuous piecewise linear
polynomials.
This element pair is widely used in computational fluid dynamics and satisfies the Ladyzhenskaya--Babuška--Brezzi
(LBB) stability condition~\cite{elman2014finite,brezzi1991}, which guarantees the stability of the discrete
problem and the absence of spurious pressure oscillations.

The implementation was carried out using the FEniCS computing platform~\cite{logg2012automated,alnaes2015},
which provides a high-level interface for automatic generation of
finite element code from the variational formulation.
To solve linear systems at each iteration, a direct method (MUMPS) is used, ensuring high accuracy for the meshes used.

\subsection{Numerical solution of the spectral problem}

Problem \eqref{eq:2-13}--\eqref{eq:2-16} corresponds to finding eigenvalues $\lambda \in \mathbb{C}$ and corresponding eigenfunctions $(\bar{\bm u}, \bar{p})$ for which
\[
\begin{split}
    \lambda \int_\Omega \bar{\bm u} \cdot \bm v \, d \bm x + 
    \frac{1}{\operatorname{Re}} 
    \int_\Omega \operatorname{grad} \bar{\bm u} : \operatorname{grad} \bm v \, d \bm x
    + \int_\Omega  \left(\operatorname{div}\bigl(\bm U \otimes \bar{\bm u}\bigr)
        + \operatorname{div}\bigl(\bar{\bm u} \otimes \bm U\bigr)\right)  \cdot \bm v \, d \bm x
    - \int_\Omega \bar{p} \, \operatorname{div}\bm v \, d\bm x & = 0 ,
    \\[2pt]
    \int_\Omega q \, \operatorname{div} \bar{\bm u} \, d \bm x &= 0 ,
\end{split}
\]
for all $(\bm v, q) \in V \times Q$.
For the finite element approximation, the same Taylor--Hood element pair
($P_2$--$P_1$) is used as for the steady problem. This ensures
consistency of the discretizations and allows, in particular, reuse of the matrices
assembled at the stage of solving the steady problem.
As a result of discretization, we obtain a generalized eigenvalue problem with a saddle-point structure of the block matrix and a degenerate mass matrix \cite{elman2014finite}.

In stability analysis, only a small number of modes determining the behavior of the flow are of interest.
We track eigenvalues with the largest real part.
To solve the partial generalized eigenvalue problem for large sparse matrices, various methods are used \cite{saad2011numerical}.
The specifics of the problems under consideration are most fully taken into account by the
Krylov--Schur algorithm~\cite{stewart2001krylov}. This method is a development of iterative projection methods on Krylov subspaces and combines
the efficiency of the Arnoldi method with a reliable restart mechanism.
The implementation was carried out using the
SLEPc (Scalable Library for Eigenvalue Problem Computations) library~\cite{hernandez2005} with the use of the FEniCS computing platform.
The SLEPc library provides efficient parallel implementations of the
Krylov--Schur algorithm and supports work with generalized eigenvalue
problems~\cite{slepc2024}.

\section{Numerical results}

The results of calculations on the study of hydrodynamic stability of two-dimensional incompressible fluid flows in a square cavity are presented.
The behavior of the dominant modes with increasing Reynolds number is studied.

\subsection{Identification of eigenvalues}

When solving the partial eigenvalue problem, it is necessary to identify the part of the spectrum of interest to us.
In the presence of imaginary eigenvalues, we have various possibilities.
Standard possibilities are related, in particular, to choosing a given number $K$ of eigenvalues with maximum (minimum) modulus.
We can also focus on the imaginary part of the eigenvalues, identifying, for example, the maximum/minimum $|\Im \lambda_k|, \ k = 1,2, ..., K$.
In stability analysis, eigenvalues with the largest real parts are of primary interest, since it is precisely the modes with eigenvalues $\Re \lambda_k > 0$ that are unstable.

\begin{figure}[ht]
	\centering
	\includegraphics[width=0.75\linewidth]{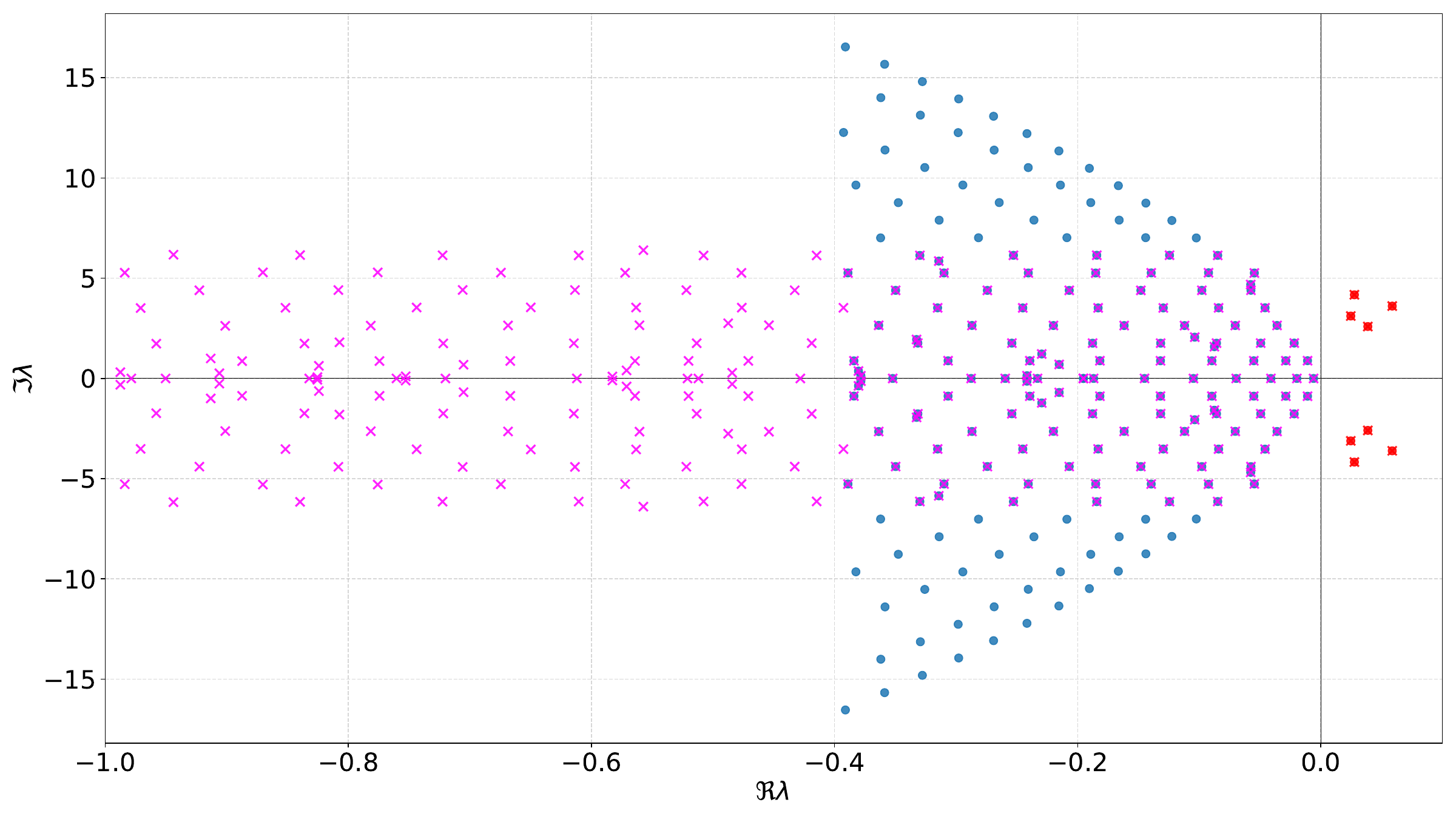}
	\caption{Eigenvalues for $\operatorname{Re} = 10\,000$ on a mesh with $N=512$: $\bullet$ --- $\max \Re \lambda_k \ (K = 200)$, $\times$ --- $\min |\lambda_k| \ (K = 1\,000)$.}
	\label{f-1}
\end{figure}

Figure~\ref{f-1} shows the eigenvalues for the problem with $\operatorname{Re} = 10\,000$ using a uniform triangular mesh with partitioning into $N=512$ segments in each direction. The part of the first $K = 1\,000$ eigenvalues with minimum modulus is marked with the symbol $\times$.
The symbol $\bullet$ marks the first $K = 200$ eigenvalues with maximum real part.
We are interested in unstable modes for which $\Re \lambda_k > 0$. In this example, we observe two pairs of unstable modes, which are highlighted in color in Fig.~\ref{f-1}. Identification of significant eigenvalues by the criterion $\max \Re \lambda_k \ (K = 100)$ was used by us in the calculations presented below.

\subsection{Mesh dependence}

\begin{figure}[ht]
	\centering
	\includegraphics[width=0.3\linewidth]{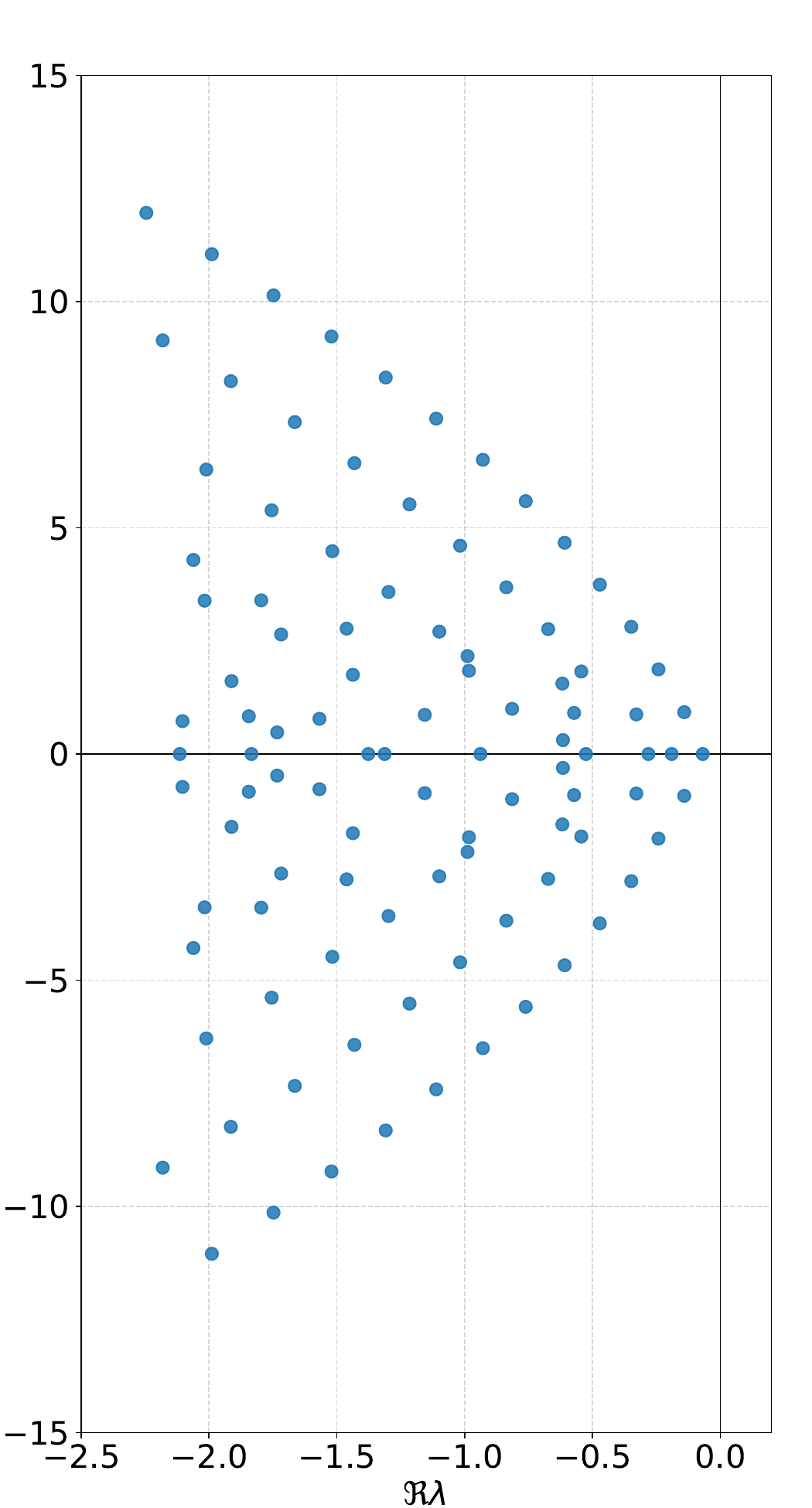}
	\includegraphics[width=0.3\linewidth]{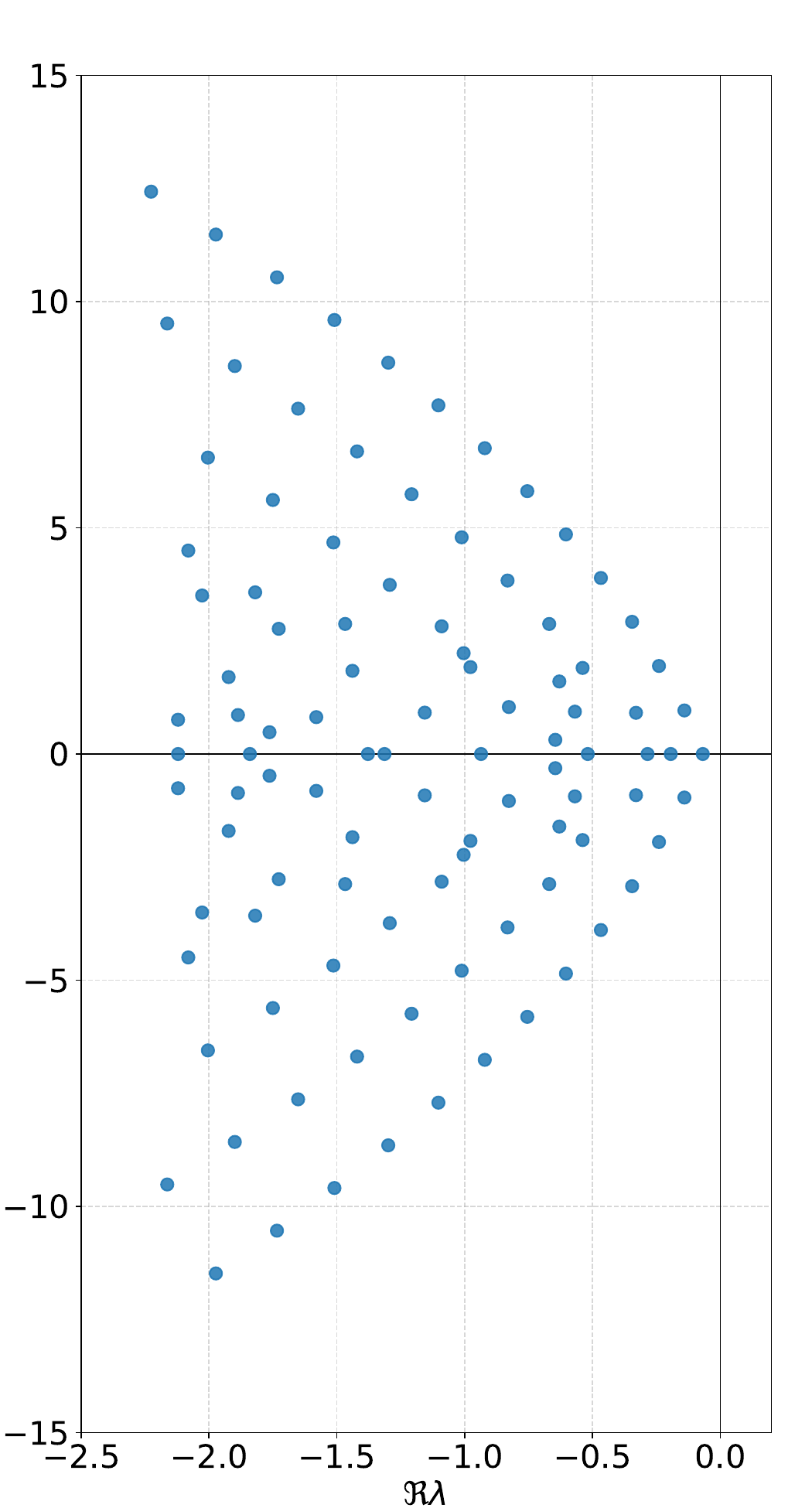}    
	\includegraphics[width=0.3\linewidth]{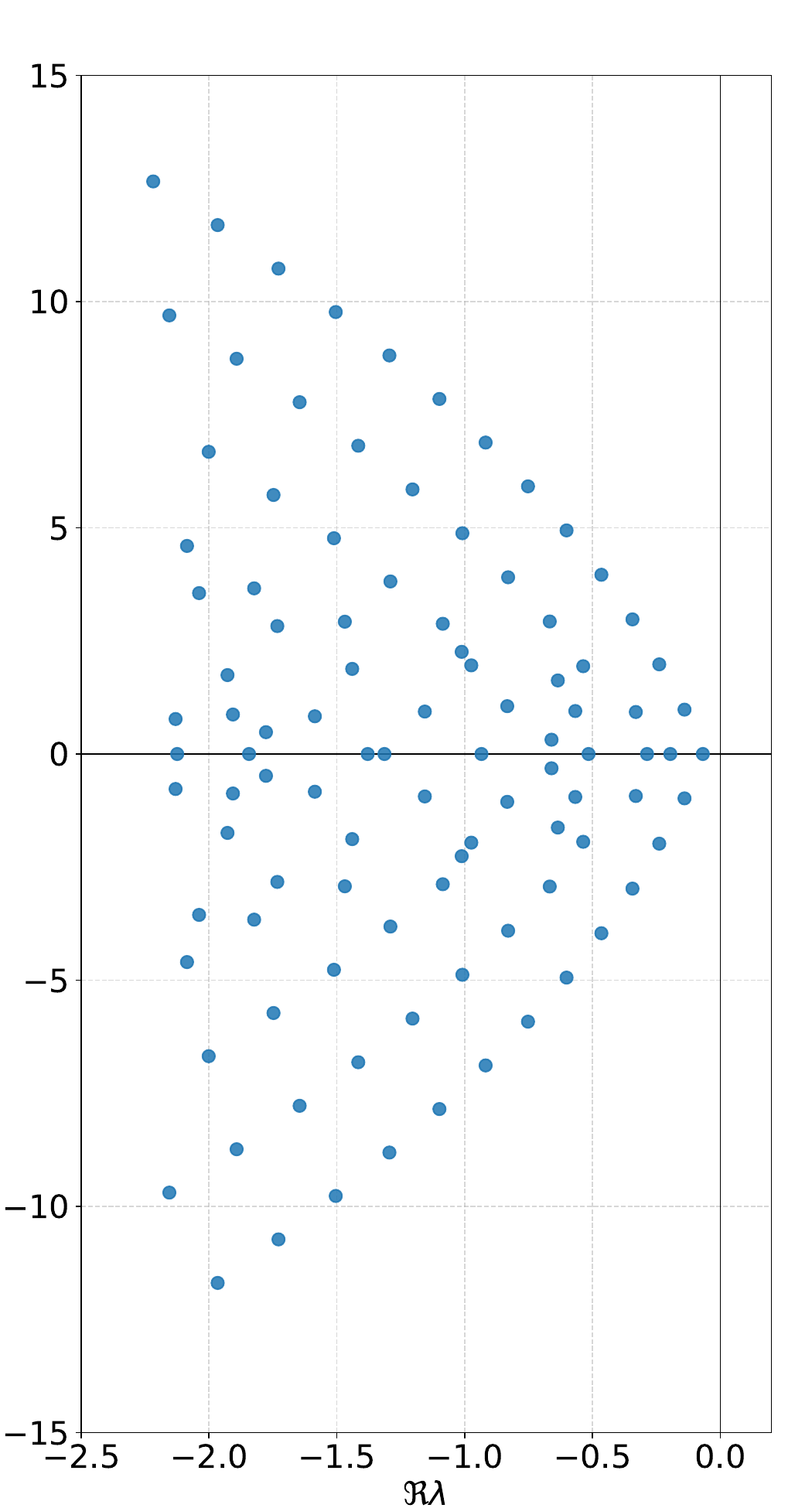}
	\caption{Eigenvalues for $\operatorname{Re} = 1\,000$ on different meshes: $N=256$ --- left, $N=512$ --- center, $N=1\,024$ --- right.}
	\label{f-2}
\end{figure}

\begin{figure}[ht]
	\centering
	\includegraphics[width=0.3\linewidth]{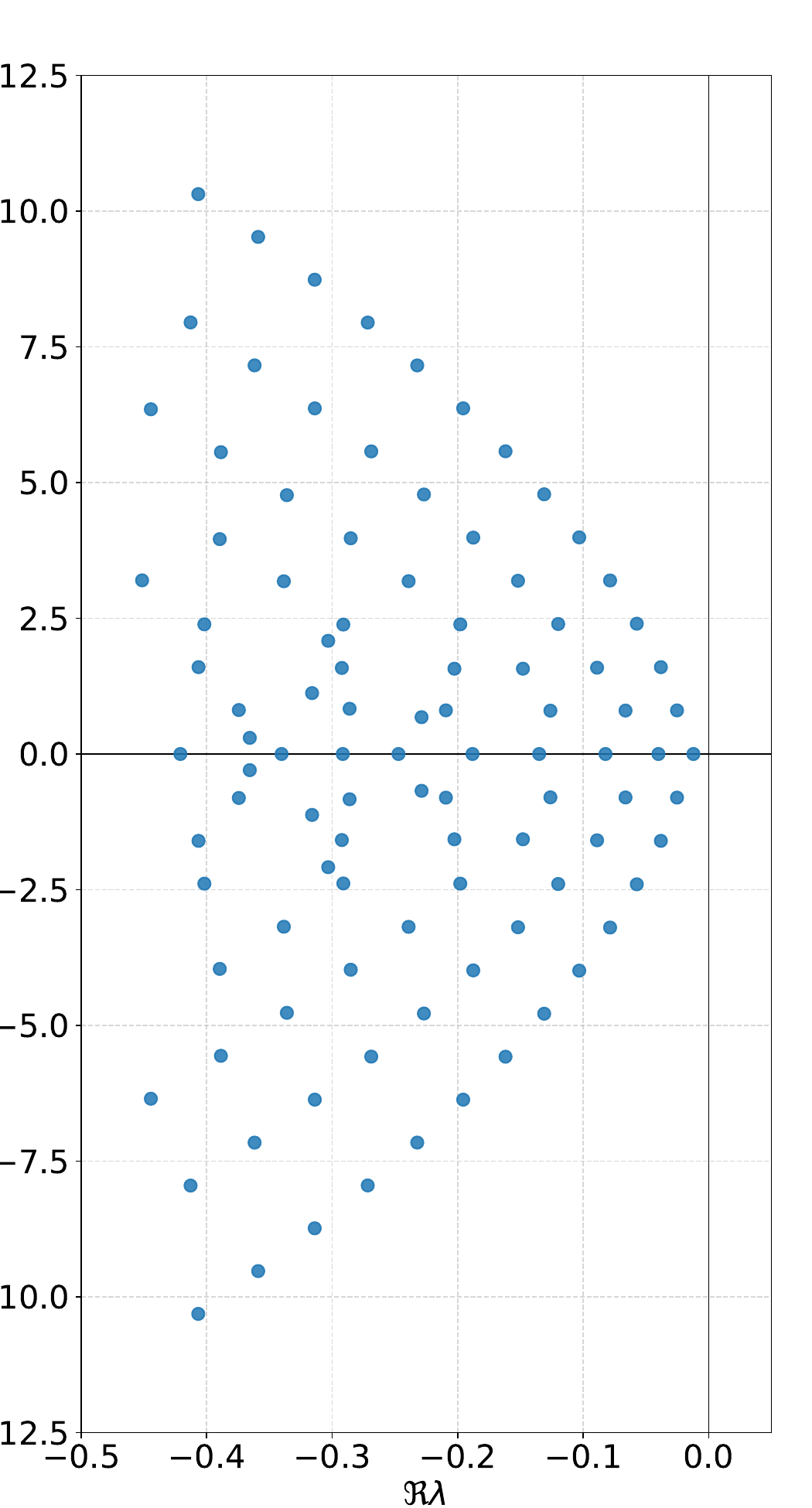}
	\includegraphics[width=0.3\linewidth]{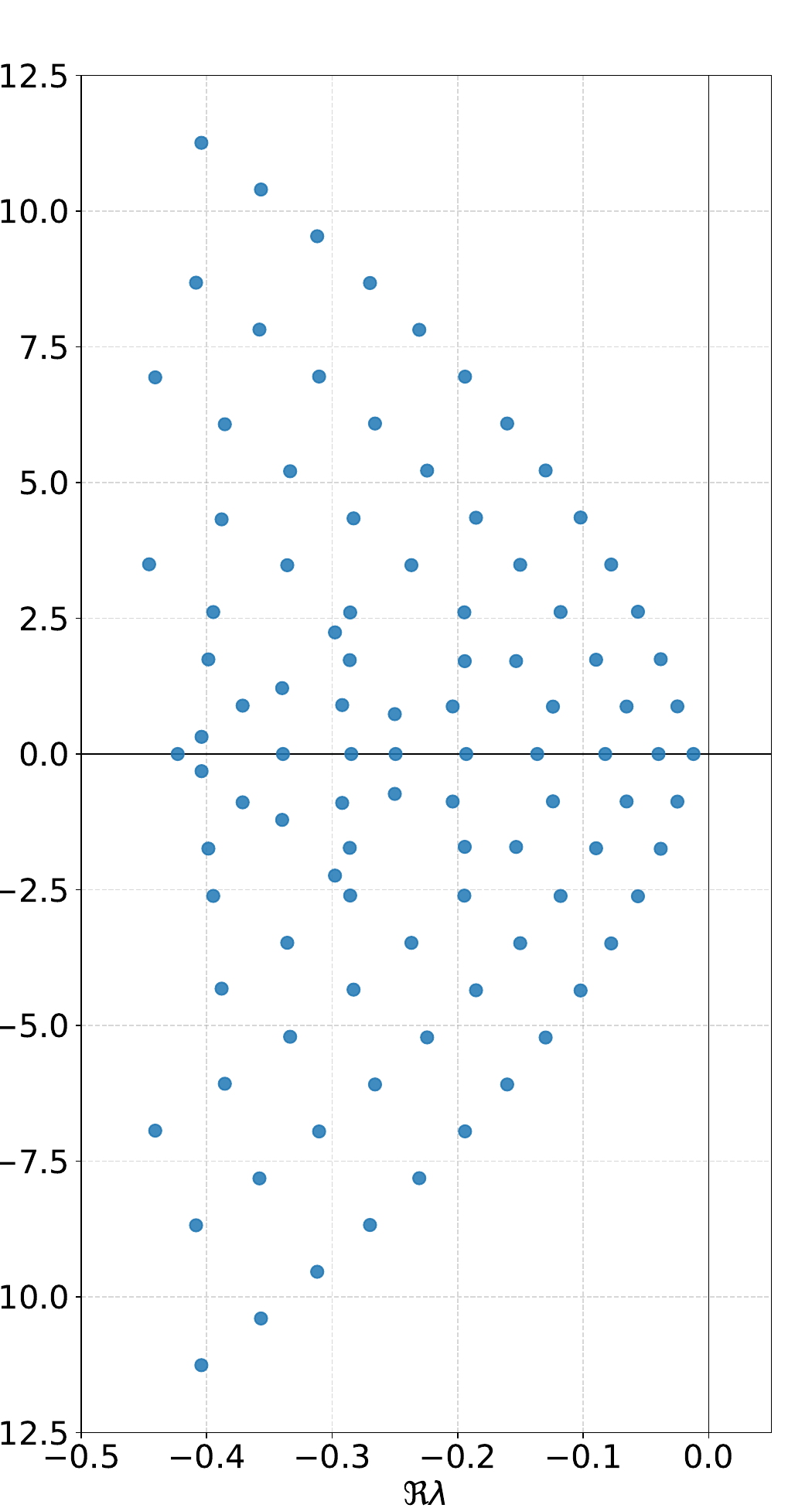}    
	\includegraphics[width=0.3\linewidth]{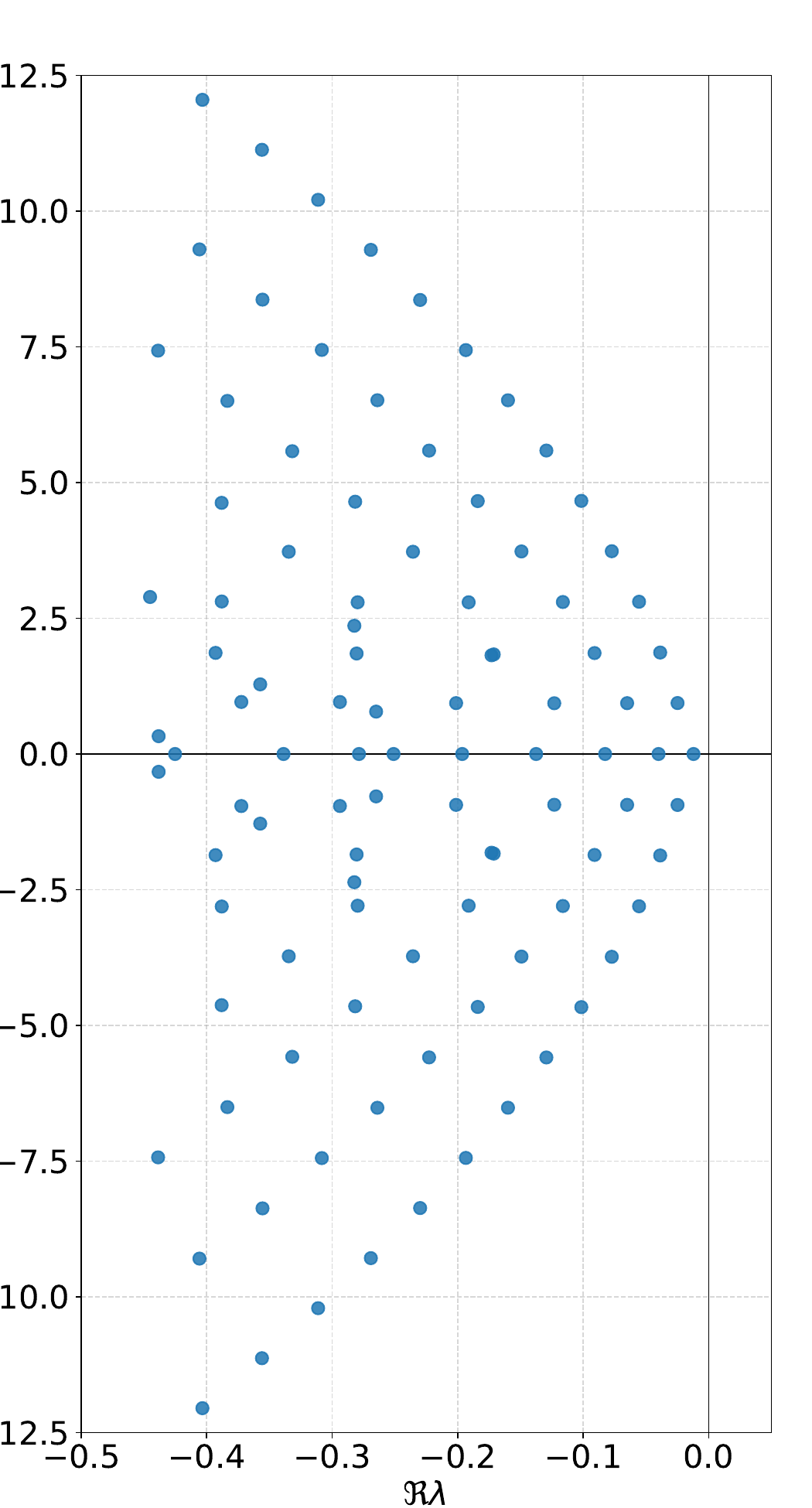}
	\caption{Eigenvalues for $\operatorname{Re} = 5\,000$ on different meshes: $N=256$ --- left, $N=512$ --- center, $N=1\,024$ --- right.}
	\label{f-3}
\end{figure}

\begin{figure}[ht]
	\centering
	\includegraphics[width=0.3\linewidth]{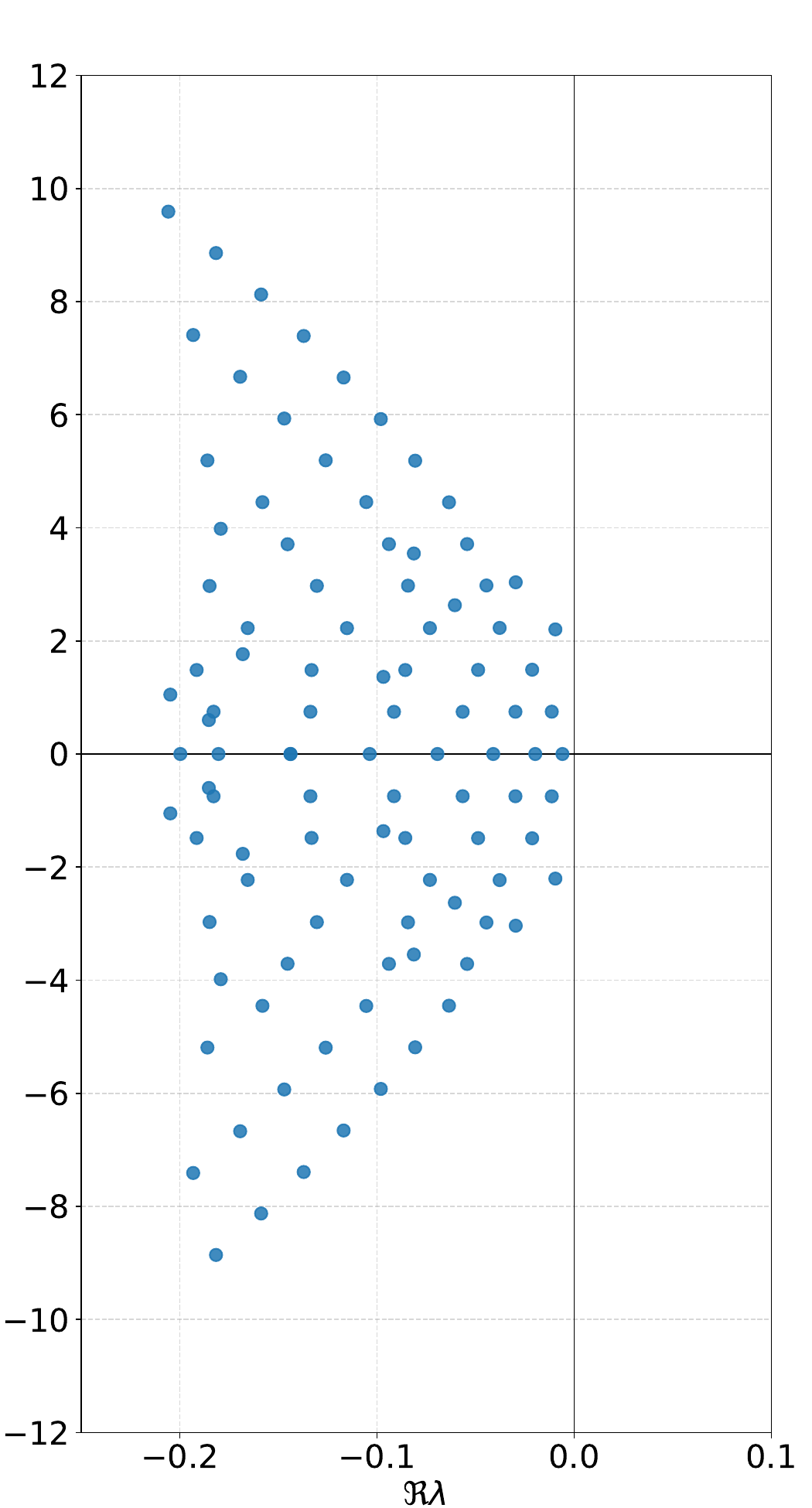}
	\includegraphics[width=0.3\linewidth]{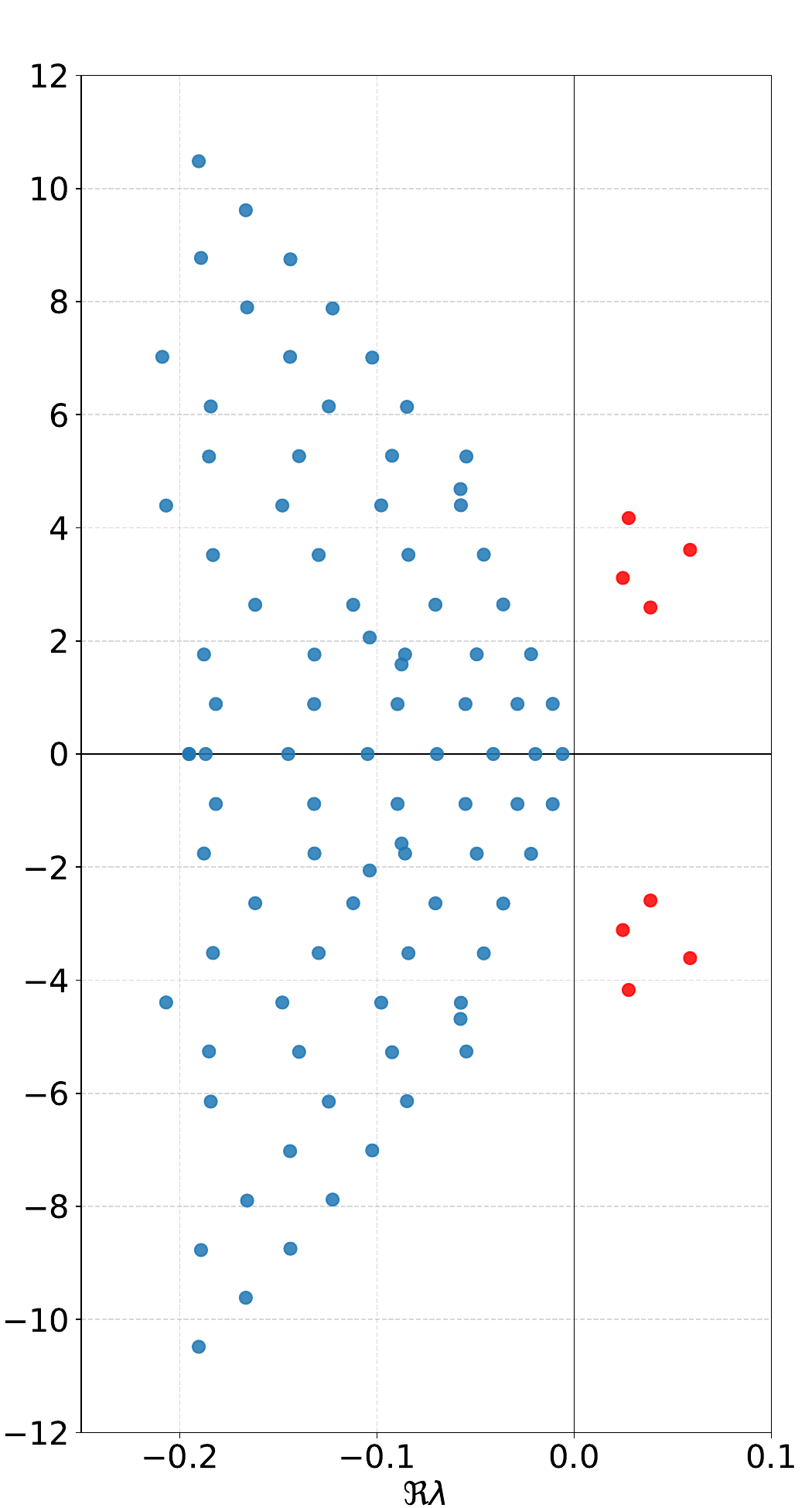}    
	\includegraphics[width=0.3\linewidth]{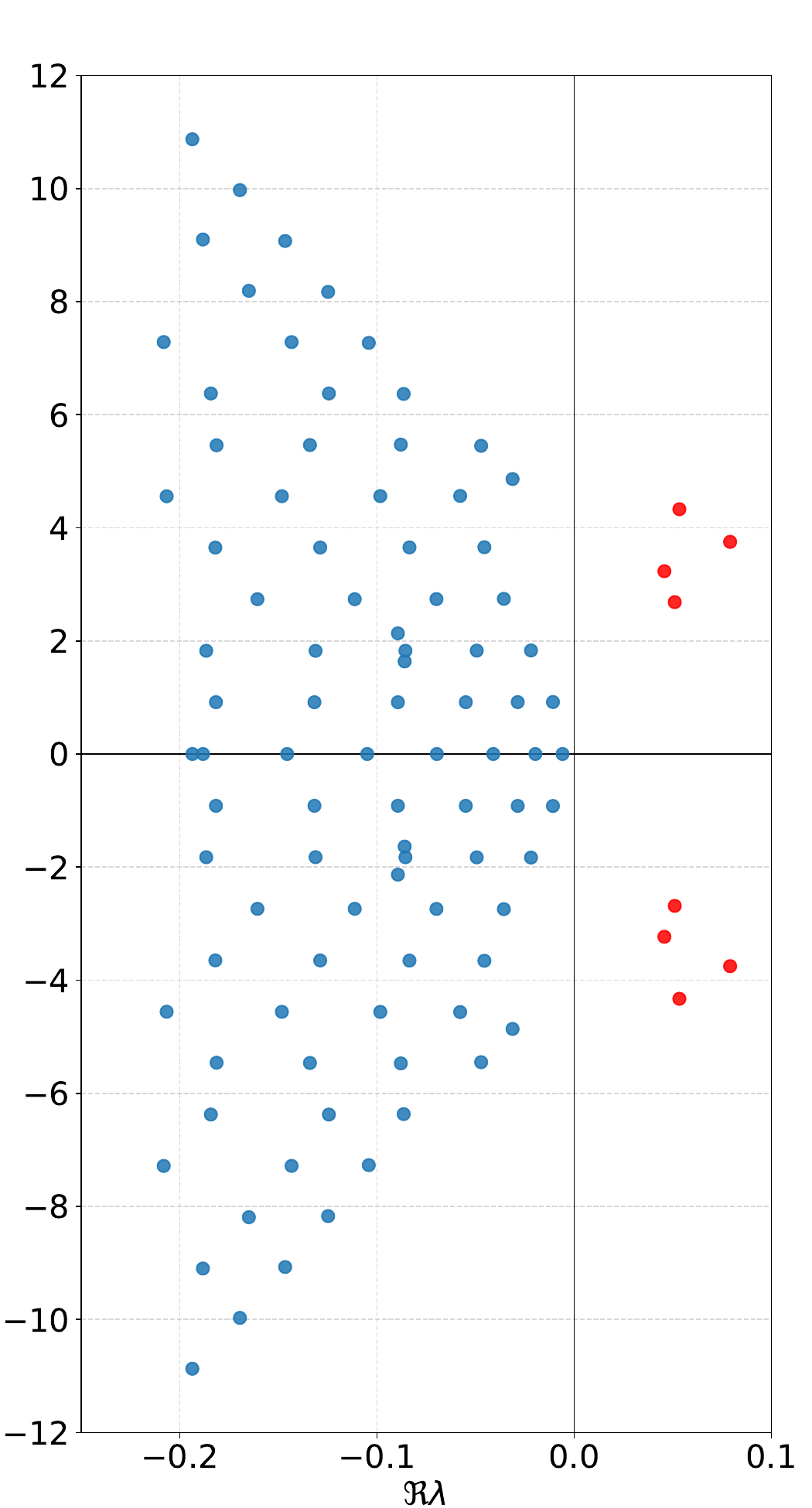}
	\caption{Eigenvalues for $\operatorname{Re} = 10\,000$ on different meshes: $N=256$ --- left, $N=512$ --- center, $N=1\,024$ --- right.}
	\label{f-4}
\end{figure}

The problem of the dependence of the dominant eigenvalues on the computational mesh is considered separately.
We performed calculations on a sequence of refined meshes from $N = 128$ to $N = 1\,024$ at various Reynolds numbers.
Figure~\ref{f-2} shows the eigenvalues for $\operatorname{Re} = 1\,000$ for three meshes: $N=256$, $N=512$, and $N=1\,024$. The first eigenvalues with maximum real part are presented in Table~\ref{t-1}.
Similar data for $\operatorname{Re} = 5\,000$ and $\operatorname{Re} = 10\,000$ are given in Fig.~\ref{f-3} and Fig.~\ref{f-4}, respectively. The numerical data are given in Table~\ref{t-2} and Table~\ref{t-3}.

\begin{table}[htbp]
	\centering
	\caption{Eigenvalues for $\operatorname{Re} = 1\,000$.}
	\label{t-1}
	\begin{tabular}{|c|l|l|l|}
		\hline
	    No. & $N = 256$ & $N = 512$ & $N = 1\,024$ \\
		\hline
		1 & $-0.06856599$ & $-0.06828967$ & $-0.06816566$ \\
		\hline
		2 & $-0.14106471 \pm 0.92460682\mathrm{i}$ & $-0.14027700 \pm 0.96243097\mathrm{i}$ & $-0.13990002 \pm 0.98068237\mathrm{i}$ \\
		\hline
		3 & $-0.19016387$ & $-0.19375443$ & $-0.19535406$ \\
		\hline
		4 & $-0.24190167 \pm 1.87101173\mathrm{i}$ & $-0.23974663 \pm 1.94589218\mathrm{i}$ & $-0.23872386 \pm 1.98200379\mathrm{i}$ \\
		\hline
		5 & $-0.28118368$ & $-0.28446167$ & $-0.28629538$ \\
		\hline
		6 & $-0.32835197 \pm 0.87497161\mathrm{i}$ & $-0.32966950 \pm 0.91087980\mathrm{i}$ & $-0.33040893 \pm 0.92837627\mathrm{i}$ \\
		\hline
		7 & $-0.34817048 \pm 2.81096120\mathrm{i}$ & $-0.34495877 \pm 2.92206602\mathrm{i}$ & $-0.34345044 \pm 2.97565537\mathrm{i}$ \\
		\hline
		8 & $-0.47111696 \pm 3.74353948\mathrm{i}$ & $-0.46697498 \pm 3.89109078\mathrm{i}$ & $-0.46504193 \pm 3.96225303\mathrm{i}$ \\
		\hline
		9 & $-0.52551873$ & $-0.51837430$ & $-0.51527402$ \\
		\hline
		10 & $-0.54361799 \pm 1.82357757\mathrm{i}$ & $-0.53853265 \pm 1.90262393\mathrm{i}$ & $-0.53617596 \pm 1.94050532\mathrm{i}$ \\
		\hline
	\end{tabular}
\end{table}

\begin{table}[htbp]
	\centering
	\caption{Eigenvalues for $\operatorname{Re} = 5\,000$.}
	\label{t-2}
	\begin{tabular}{|c|l|l|l|}
		\hline
	    No. & $N = 256$ & $N = 512$ & $N = 1\,024$ \\
		\hline
		1 & $-0.01214438$ & $-0.01210141$ & $-0.01206677$ \\
		\hline
		2 & $-0.02520852 \pm 0.80299420\mathrm{i}$ & $-0.02490945 \pm 0.87691695\mathrm{i}$ & $-0.02480186 \pm 0.93881034\mathrm{i}$ \\
		\hline
		3 & $-0.03806527 \pm 1.59957072\mathrm{i}$ & $-0.03819202 \pm 1.74612555\mathrm{i}$ & $-0.03861312 \pm 1.86894004\mathrm{i}$ \\
		\hline
		4 & $-0.04002541$ & $-0.03997646$ & $-0.03992780$ \\
		\hline
		5 & $-0.05726808 \pm 2.40005217\mathrm{i}$ & $-0.05633290 \pm 2.62035575\mathrm{i}$ & $-0.05545044 \pm 2.80487079\mathrm{i}$ \\
		\hline
		6 & $-0.06625268 \pm 0.80031720\mathrm{i}$ & $-0.06539835 \pm 0.87438623\mathrm{i}$ & $-0.06496877 \pm 0.93640832\mathrm{i}$ \\
		\hline
		7 & $-0.07853635 \pm 3.19620836\mathrm{i}$ & $-0.07765282 \pm 3.48943466\mathrm{i}$ & $-0.07719464 \pm 3.73517197\mathrm{i}$ \\
		\hline
		8 & $-0.08219743$ & $-0.08240601$ & $-0.08251134$ \\
		\hline
		9 & $-0.08891307 \pm 1.58975193\mathrm{i}$ & $-0.08969499 \pm 1.73624905\mathrm{i}$ & $-0.09092916 \pm 1.85915921\mathrm{i}$ \\
		\hline
		10 & $-0.10313377 \pm 3.99048175\mathrm{i}$ & $-0.10204817 \pm 4.35628558\mathrm{i}$ & $-0.10145548 \pm 4.66273295\mathrm{i}$ \\
		\hline
	\end{tabular}
\end{table}

\begin{table}[htbp]
	\centering
	\caption{Eigenvalues for $\operatorname{Re} = 10\,000$.}
	\label{t-3}
	\begin{tabular}{|c|l|l|l|}
		\hline
	    No. & $N = 256$ & $N = 512$ & $N = 1\,024$ \\
		\hline
		1 & $-0.00595383$ & $\hphantom{X} 0.05876709 \pm 3.60987961\mathrm{i}$ & $\hphantom{X} 0.07906966 \pm 3.75144934\mathrm{i}$ \\
		\hline
		2 & $-0.00954666 \pm 2.20336105\mathrm{i}$ & $\hphantom{X} 0.03867942 \pm 2.59064987\mathrm{i}$ & $\hphantom{X} 0.05334142 \pm 4.32885208\mathrm{i}$ \\
		\hline
		3 & $-0.01137723 \pm 0.74846661\mathrm{i}$ & $\hphantom{X} 0.02763670 \pm 4.17235153\mathrm{i}$ & $\hphantom{X} 0.05098806 \pm 2.68558842\mathrm{i}$ \\
		\hline
		4 & $-0.01974827$ & $\hphantom{X} 0.02467004 \pm 3.11335945\mathrm{i}$ & $\hphantom{-} 0.04572733 \pm 3.23200591\mathrm{i}$ \\
		\hline
		5 & $-0.02131235 \pm 1.49054876\mathrm{i}$ & $-0.00592648$ & $-0.00592208$ \\
		\hline
		6 & $-0.02963922 \pm 3.03640907\mathrm{i}$ & $-0.01089822 \pm 0.88634546\mathrm{i}$ & $-0.01078486 \pm 0.91929850\mathrm{i}$ \\
		\hline
		7 & $-0.02979813 \pm 0.74714104\mathrm{i}$ & $-0.01969248$ & $-0.01968453$ \\
		\hline
		8 & $-0.03780683 \pm 2.22966823\mathrm{i}$ & $-0.02186675 \pm 1.76516542\mathrm{i}$ & $-0.02189836 \pm 1.83094244\mathrm{i}$ \\
		\hline
		9 & $-0.04104718$ & $-0.02877523 \pm 0.88465942\mathrm{i}$ & $-0.02864613 \pm 0.91751676\mathrm{i}$ \\
		\hline
		10 & $-0.04449695 \pm 2.98002761\mathrm{i}$ & $-0.03598147 \pm 2.64424484\mathrm{i}$ & $-0.03124785 \pm 4.86200511\mathrm{i}$ \\
		\hline
	\end{tabular}
\end{table}

For $\operatorname{Re} = 1\,000$, the leading eigenvalues converge well as the mesh is refined. Moreover, all modes are stable ($\Re \lambda_k < 0$).
For $\operatorname{Re} = 5\,000$, a significant influence of the mesh on the imaginary parts of the eigenvalues is noted.
For $\operatorname{Re} = 10\,000$ on a coarse mesh ($N=256$), there are no unstable modes, while for $N=512$ and $N=1\,024$, 4 pairs of unstable modes appear.

Thus, the calculated values for the dominant eigenvalues depend quantitatively and qualitatively on the mesh size. As the Reynolds number increases, the decrease in accuracy manifests itself more strongly.
For an acceptable quantitative estimate, sufficiently detailed computational meshes must be used.

\subsection{Dependence on the Reynolds number}

The analysis of how the spectral characteristics of the flow depend on the Reynolds number is particularly important. For this purpose, a series of parametric calculations was performed on a mesh with $N = 1\,024$ in the range $\operatorname{Re} \in [8\,000, 10\,500]$. For each value of the Reynolds number, the eigenvalue problem was solved and the modes responsible for loss of stability were determined.

The results agree with published estimates of the critical Reynolds number for the two-dimensional cavity problem. According to the literature sources noted above, the first Hopf bifurcation for the two-dimensional steady solution occurs at $\operatorname{Re}_{\mathrm{cr}} \approx 8\,000$.
As the Reynolds number increases, a successive increase in the number of unstable modes is observed.

\begin{figure}[htbp]
	\centering
	\begin{subfigure}[b]{0.32\linewidth}
		\centering
		\includegraphics[width=\linewidth]{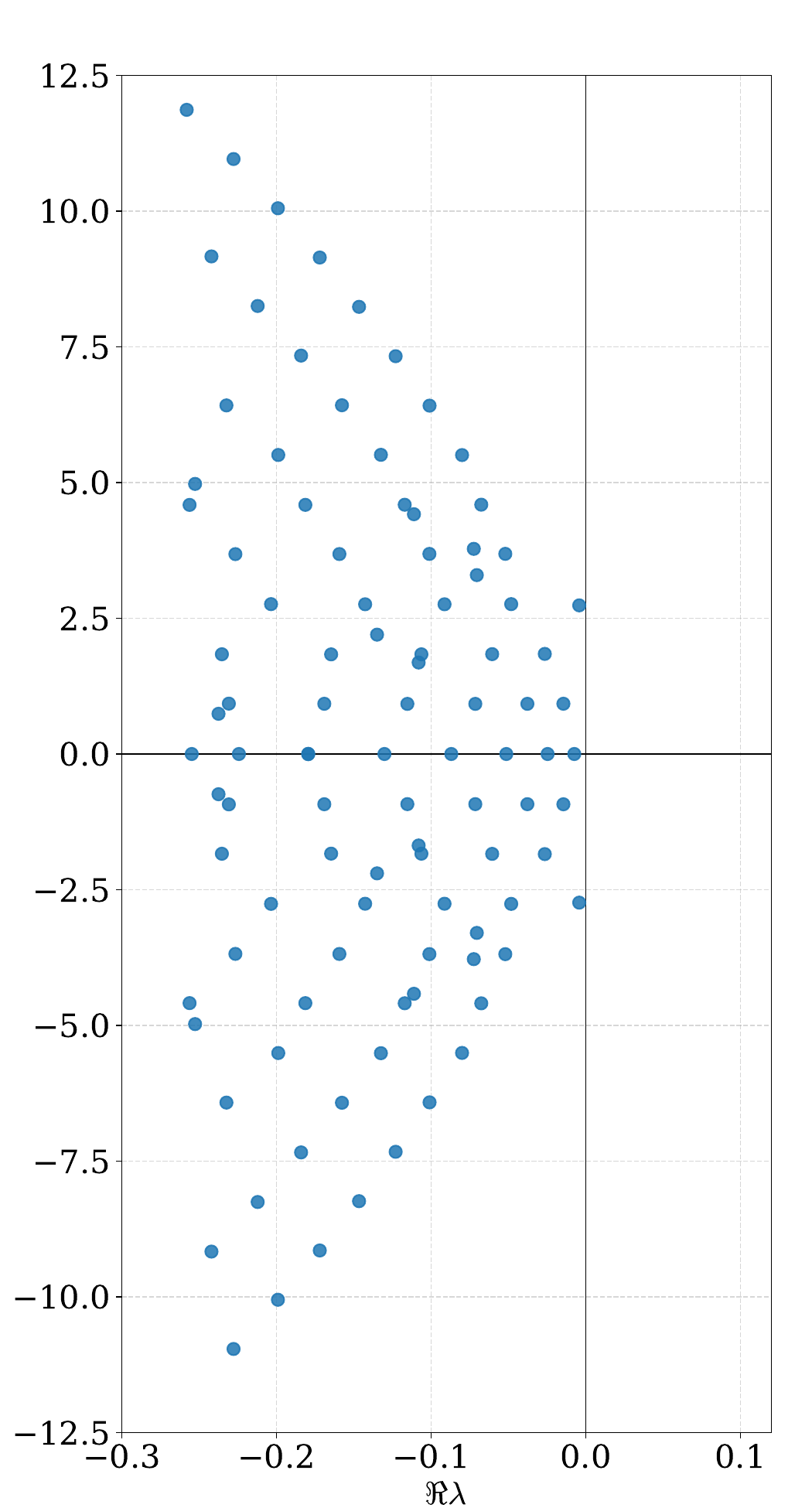}
		\caption{}
		\label{f-5-a}
	\end{subfigure}\hfill
	\begin{subfigure}[b]{0.32\linewidth}
		\centering
		\includegraphics[width=\linewidth]{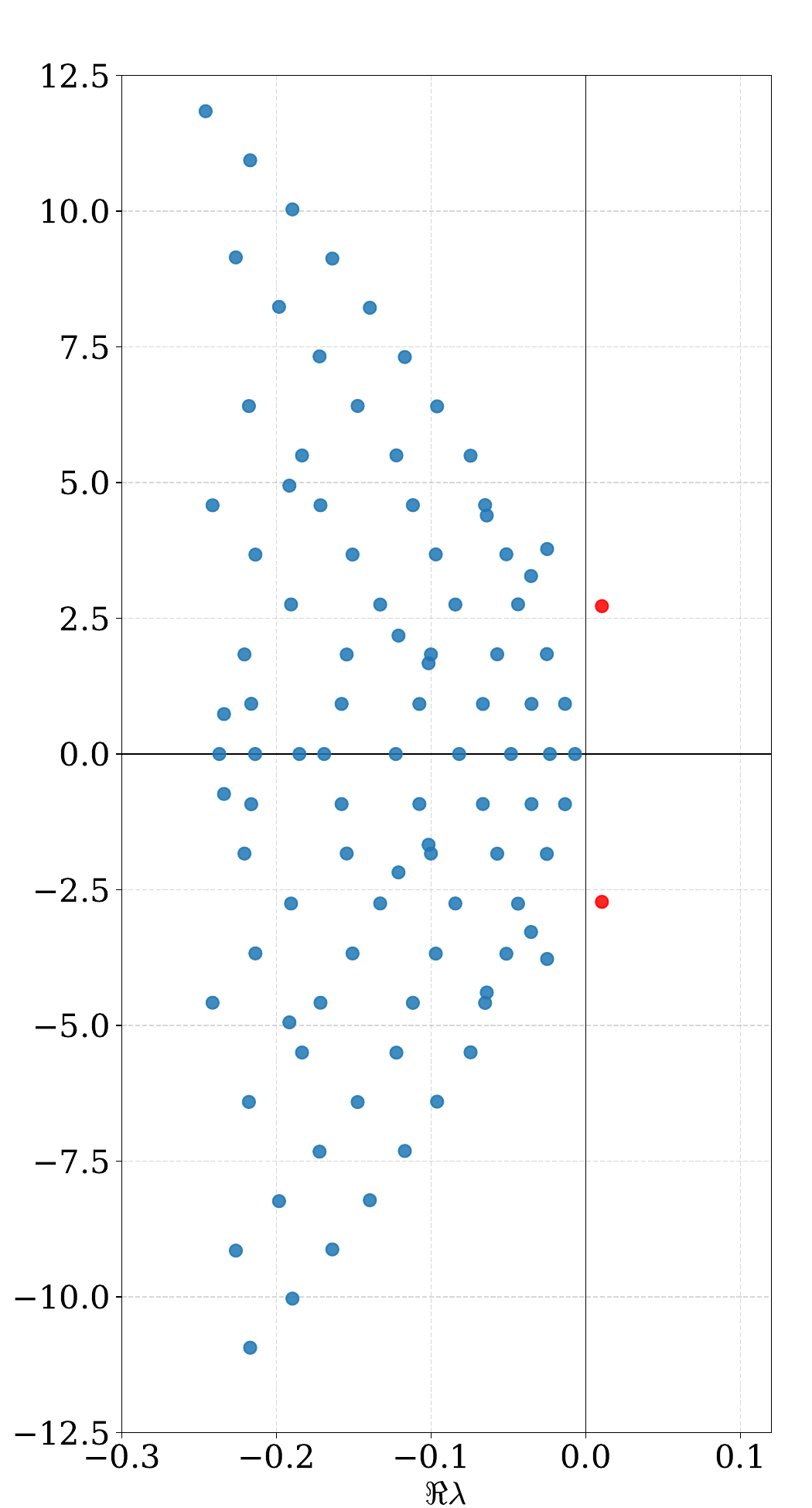}
		\caption{}
		\label{f-5-b}
	\end{subfigure}\hfill
	\begin{subfigure}[b]{0.32\linewidth}
		\centering
		\includegraphics[width=\linewidth]{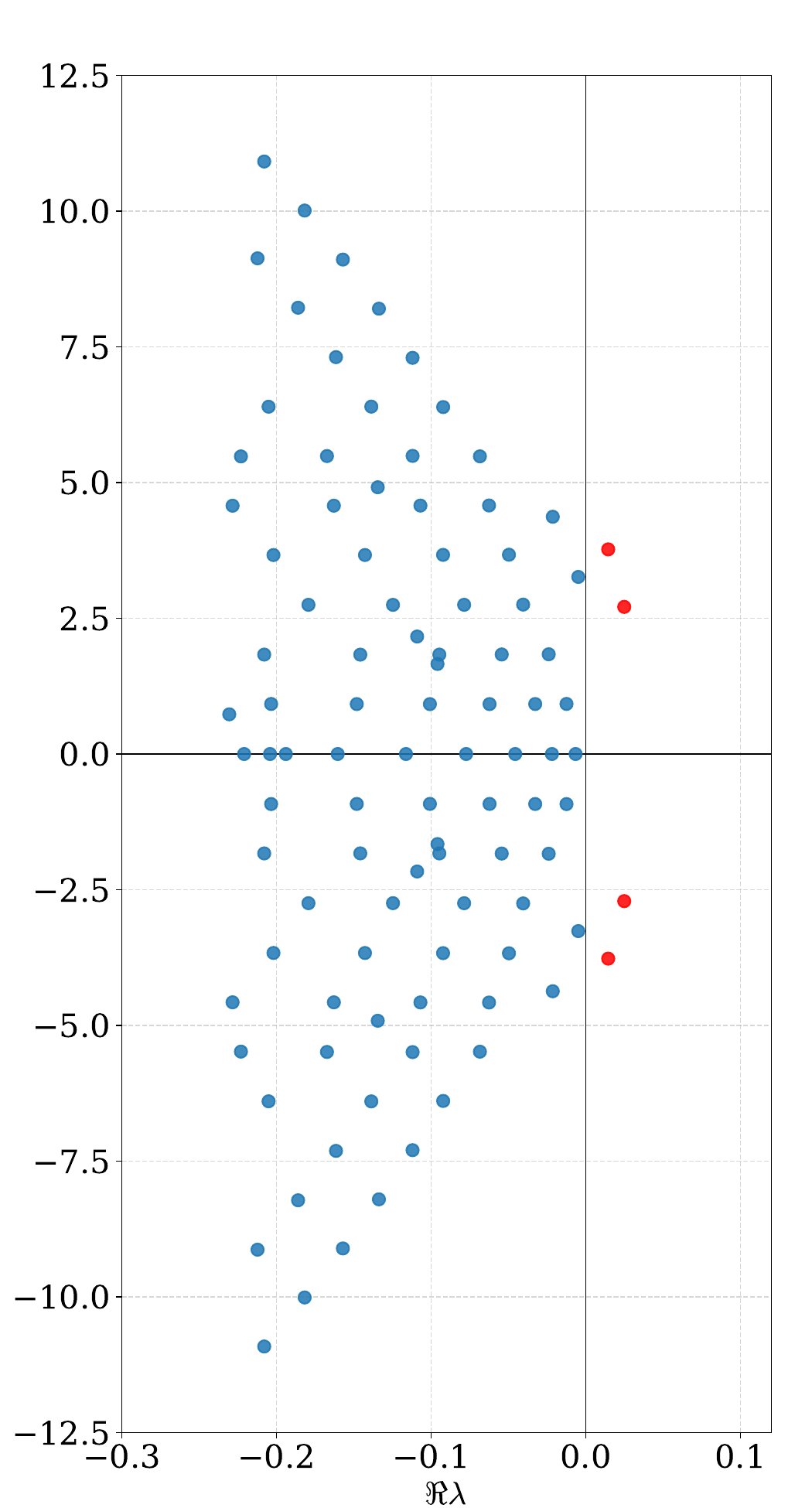}
		\caption{}
		\label{f-5-c}
	\end{subfigure}
	
	\vspace{0.5em}
	
	\begin{subfigure}[b]{0.32\linewidth}
		\centering
		\includegraphics[width=\linewidth]{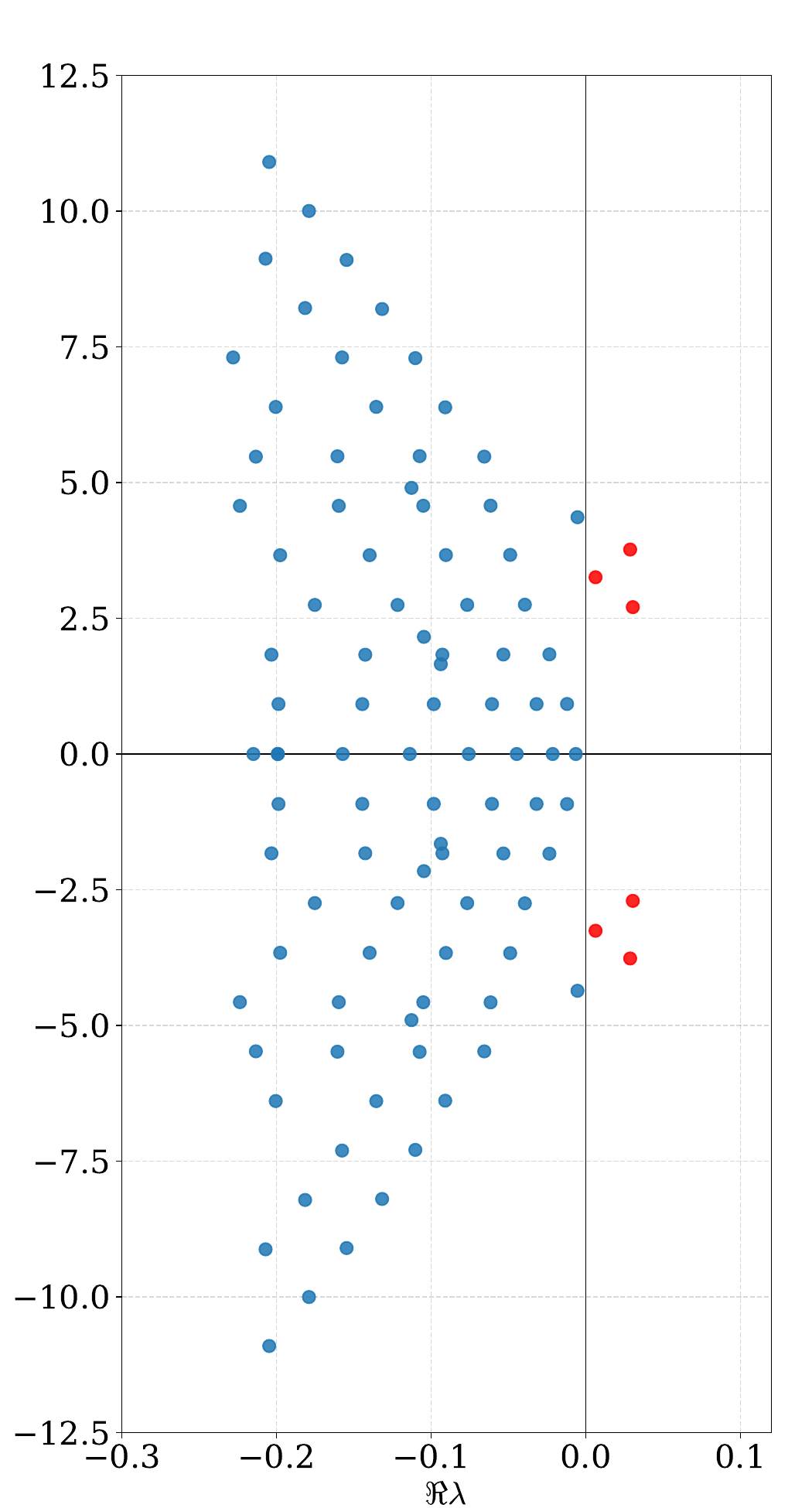}
		\caption{}
		\label{f-5-d}
	\end{subfigure}\hfill
	\begin{subfigure}[b]{0.32\linewidth}
		\centering
		\includegraphics[width=\linewidth]{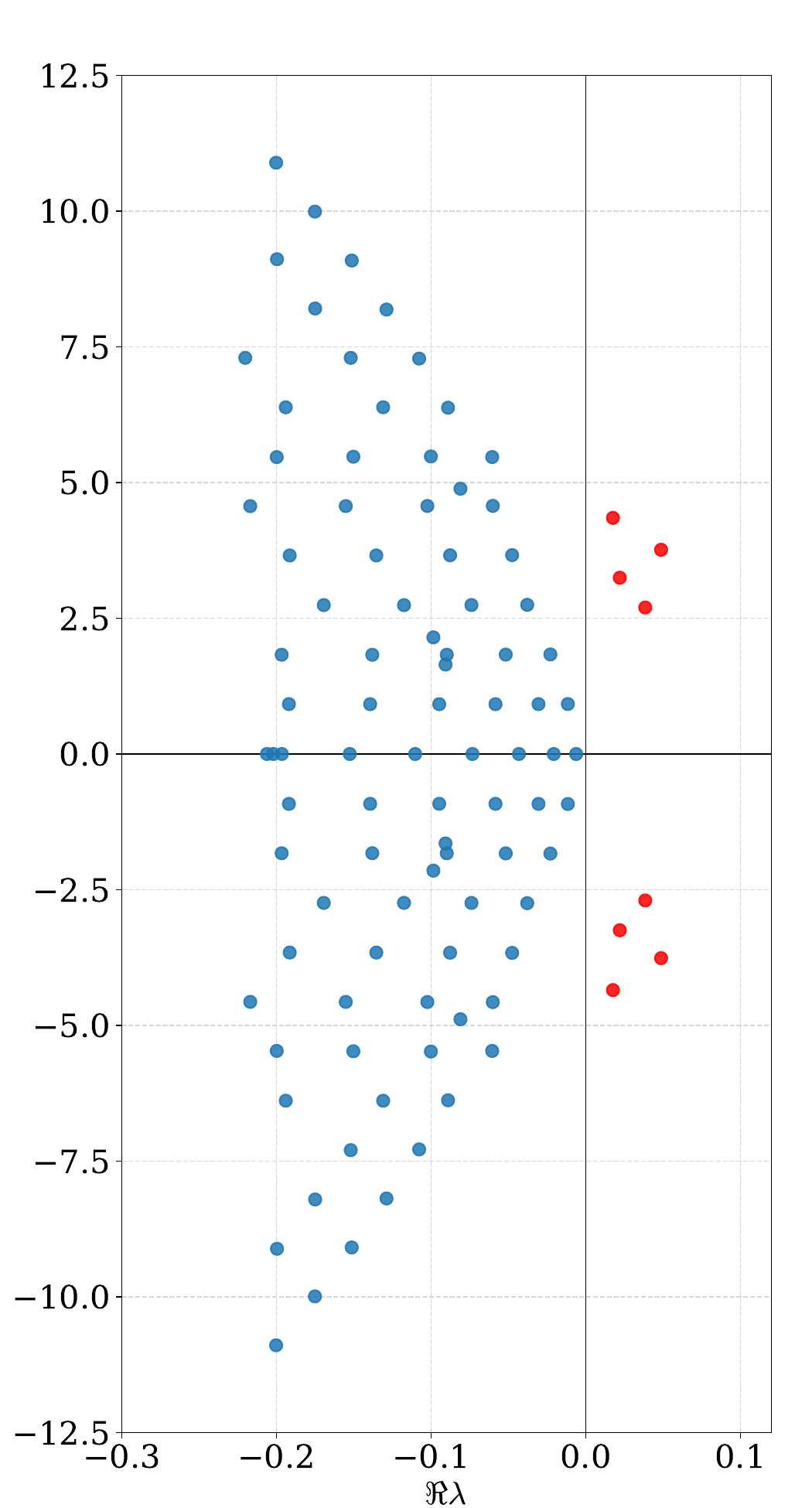}
		\caption{}
		\label{f-5-e}
	\end{subfigure}\hfill
	\begin{subfigure}[b]{0.32\linewidth}
		\centering
		\includegraphics[width=\linewidth]{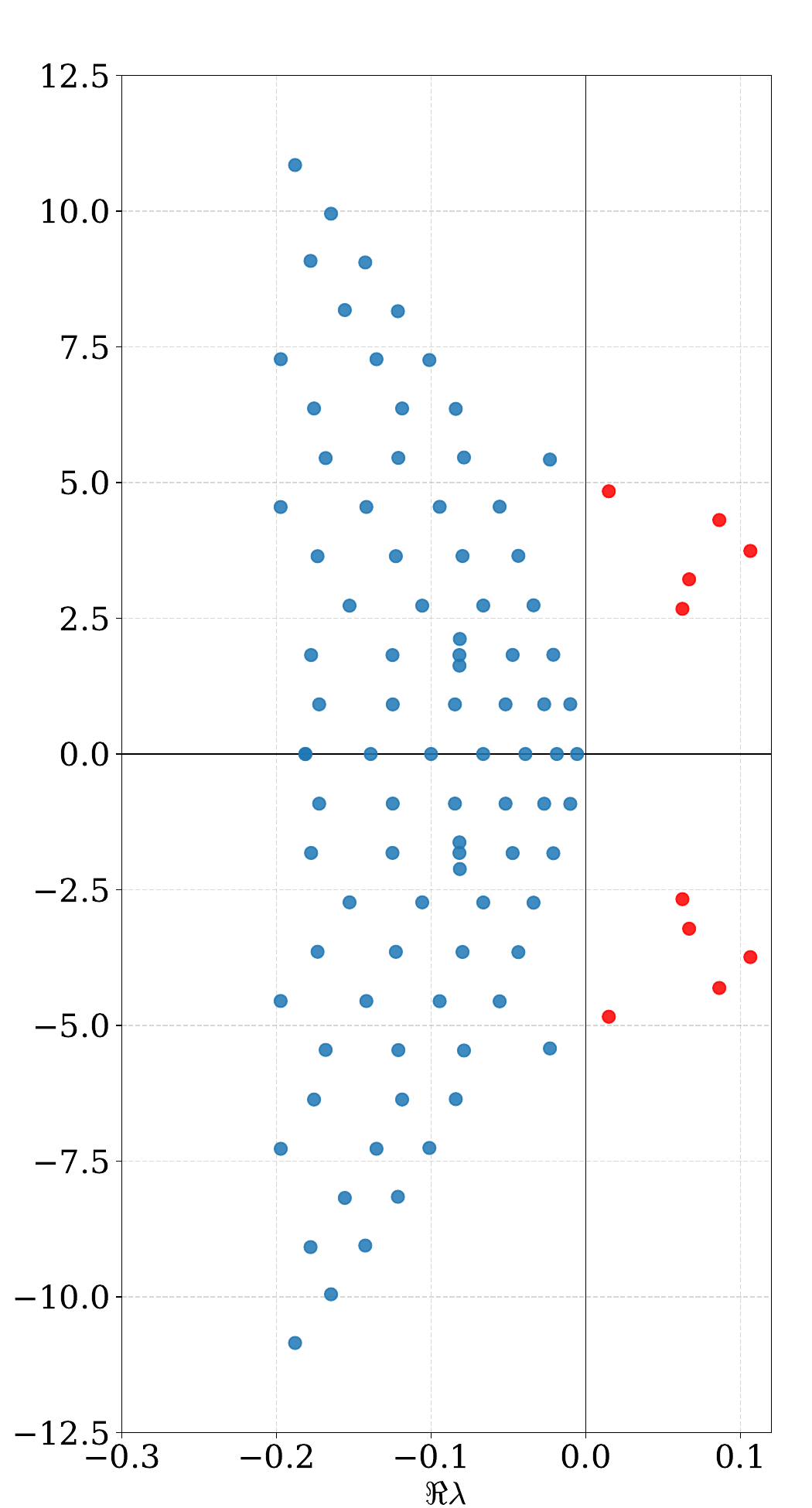}
		\caption{}
		\label{f-5-f}
	\end{subfigure}
	
	\caption{Eigenvalues on a mesh $N=1\,024$ at various Reynolds numbers: (a)~--- $\operatorname{Re} = 8\,000$; (b)~--- $\operatorname{Re} = 8\,500$; (c)~--- $\operatorname{Re} = 9\,000$; (d)~--- $\operatorname{Re} = 9\,200$; (e)~--- $\operatorname{Re} = 9\,500$;
	(f)~---$\operatorname{Re} = 10\,500$.}
	\label{f-5}
\end{figure}

Figure~\ref{f-5} shows the dominant modes for several representative values of $\operatorname{Re}$. The qualitative evolution of the spectrum is as follows.
At $\operatorname{Re} = 8\,000$, all computed eigenvalues lie in the left half-plane ($\Re \lambda < 0$), which means stability of the steady flow with respect to two-dimensional perturbations.
At $\operatorname{Re} = 8\,500$, one pair of complex-conjugate eigenvalues with positive real part appears in the spectrum. This corresponds to a Hopf bifurcation and the occurrence of a periodic regime. The imaginary part of the critical pair determines the oscillation frequency after the bifurcation.
With a further increase in the Reynolds number, the number of unstable mode pairs increases. At $\operatorname{Re} = 10\,500$, five pairs of unstable modes are already present in the spectrum, which indicates a sequence of secondary bifurcations and increasing complexity in the flow dynamics.

\begin{figure}[htbp]
	\centering
	\includegraphics[width=0.75\linewidth]{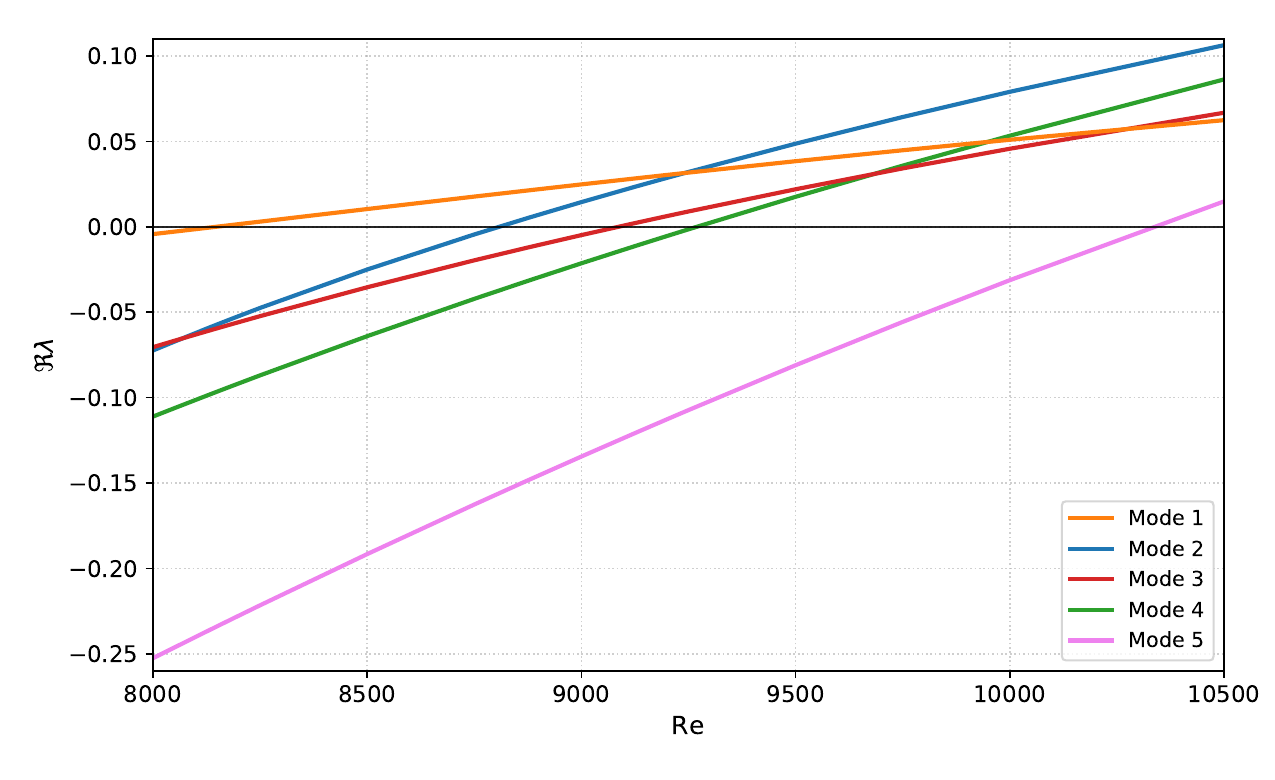}
	\includegraphics[width=0.75\linewidth]{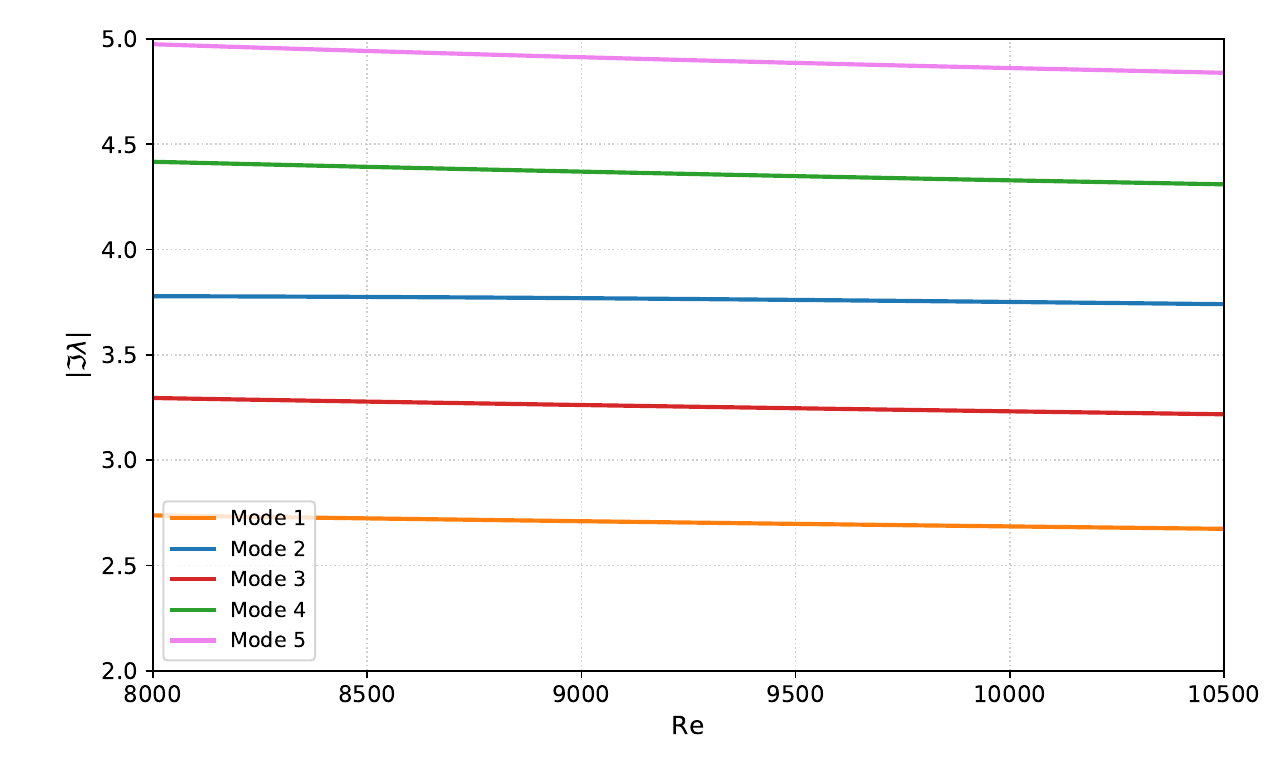}
	\caption{Dependence of the real and imaginary parts of the leading eigenvalues on $\operatorname{Re}$.}
	\label{f-6}
\end{figure}

The main result of our study is an analysis of how the spectral characteristics of the dominant modes depend on the Reynolds number. The corresponding data for the real and imaginary parts of these eigenvalues are shown in Fig.~\ref{f-6}. These calculated data allow one to estimate the critical Reynolds numbers at which the stability changes. The estimated critical Reynolds numbers for the first five modes are summarized in Table~\ref{t-4}.

\clearpage

\begin{table}[htbp]
	\centering
	\caption{Critical Reynolds numbers.}
	\label{t-4}
	\begin{tabular}{|c|c|c|c|c|c|}
		\hline
	    No. & $1$ & $2$ & $3$  & $4$ & $5$ \\
		\hline
		$\operatorname{Re}_{\mathrm{cr}}$ &  $8\,146$ &  $8\,807$ & $9\,086$   & $9\,269$  & $10\,338$ \\
		\hline
	\end{tabular}
\end{table}

It should be emphasized that the specific numerical values of the critical Reynolds numbers, as was shown earlier, depend on the computational mesh and therefore are subject to refinement in more accurate calculations. In particular, in the literature the spread of estimates even for the first bifurcation value remains significant. Nevertheless, the qualitative behavior of the spectrum itself, associated with a successive increase in the number of unstable modes with increasing Reynolds number, agrees well with expectations.

\section{Conclusion}

\begin{enumerate}[(1)]

\item This work presents a numerical study of the linear hydrodynamic
stability of steady flows of a viscous incompressible fluid in a square
cavity with a moving upper lid. Particular attention is paid to the systematic
analysis of the evolution of the spectrum of the linearized operator with a successive
increase in the Reynolds number up to $\operatorname{Re} = 10\,500$.

\item Steady solutions were constructed by the Picard iterative method with linearization of the convective term of a special form and relaxation.
Such an algorithm ensures reliable convergence at high Reynolds numbers and allows the steady branch to be continued into the instability region. Spatial discretization was performed by the finite element method on triangular meshes with
the stable Taylor--Hood pair. The spectral
problem was solved by the Krylov--Schur algorithm (SLEPc library) and the eigenvalues with the largest real parts were identified.

\item The dependence of the determining eigenvalues on the computational mesh was studied.
It is shown that as the Reynolds number increases, the sensitivity of the spectrum to spatial resolution increases.
The main parametric calculations were performed on a uniform triangular mesh with $N = 1\,024$ subdivisions in each coordinate direction.

\item The obtained results qualitatively agree with known data on the first loss of stability.
A sequence of secondary bifurcations is observed.
At $\operatorname{Re} = 8\,000$, no unstable modes are found;
at $\operatorname{Re} = 8\,500$, there is one pair of
complex-conjugate eigenvalues with positive real part (Hopf bifurcation), and at $\operatorname{Re} = 10\,500$, the number of unstable
mode pairs reaches five.

\item The main result of the work is the construction, in a unified spectral formulation, of a picture of the evolution of the dominant eigenmodes of the steady solution with increasing Reynolds number and the identification of the modes responsible for successive losses of stability.
The obtained critical values refer to two-dimensional perturbations. Extension of the approach to the three-dimensional case and refinement of the critical numbers on refined meshes represent natural directions for further research.

\end{enumerate}

\end{document}